\documentclass[12pt,a4paper,reqno]{amsart}
\usepackage[utf8]{inputenc}
\usepackage[T1]{fontenc}
\usepackage[foot]{amsaddr}
\usepackage{amsmath}
\usepackage{amsfonts}
\usepackage{amssymb}
\usepackage{algorithm}
\usepackage{algpseudocode}
\usepackage{xcolor}
\usepackage{soul}
\usepackage{graphicx}
\usepackage{wrapfig}
\usepackage{subcaption}
\usepackage[nopar]{lipsum}
\usepackage{amsthm}
\usepackage{amssymb}
\usepackage{euscript}
\usepackage{blkarray}
\usepackage{tikz-cd}
\usepackage{makecell, multirow, tabularx,booktabs}
\usepackage{enumitem}
\setlist{nolistsep}
\usepackage{mathtools}
\usepackage{soul}
\usepackage{amscd}
\usepackage{float}
\usepackage[pagebackref,colorlinks=true,
linkcolor=blue,
urlcolor=red,colorlinks,
citecolor=green]{hyperref}
\usepackage{cleveref}
\usepackage{geometry}
\usepackage{diagbox}
\usepackage{pgf,tikz,pgfplots}
\tikzstyle{none}=[]
\tikzstyle{new style 0}=[draw,circle,fill=white]
\tikzstyle{new edge style 1}=[draw,dashed]
\tikzstyle{new edge style 1}=[draw,dashed]
\usepackage{adjustbox}
\usepackage{setspace}
\usepackage[numbers,sort&compress]{natbib}

\allowdisplaybreaks[4]
\usetikzlibrary{calc}
\usetikzlibrary{patterns}
\pgfplotsset{compat=1.14}
\numberwithin{equation}{section}
\newtheorem{thm}{Theorem}[section]
\newtheorem{prop}[thm]{Proposition}

\newtheorem{cor}[thm]{Corollary}
\newtheorem{exm}[thm]{Example}
\newtheorem{df}[thm]{Definition}
\newtheorem{rem}[thm]{Remark}
\newcommand {\R} {\mathbb{R}}

\newcommand{\mf}{\mathfrak}

\makeatletter
\@namedef{subjclassname@1991}{Subject}
\@namedef{subjclassname@2000}{Subject}
\@namedef{subjclassname@2010}{Subject}
\@namedef{subjclassname@2020}{Subject}
\makeatother
\begin{document}
	\title[Hypergraph Symmetries and Invariant Subspaces of Related Matrices ]{
 Symmetries of Hypergraphs and Some Invariant Subspaces of Matrices Associated with Hypergraphs}

	\author[Banerjee]{Anirban Banerjee} 
	\email[Banerjee ]{\textit {{\scriptsize anirban.banerjee@iiserkol.ac.in}}}
 \address[Banerjee]{Department of Mathematics and Statistics, Indian Institute of Science Education and Research Kolkata, Mohanpur-741246, India}
		\author[ Parui]{Samiron Parui} 
	\email[ Parui ]{\textit {{\scriptsize  samironparui@gmail.com}}}
\address[Parui]{School of Mathematical Sciences,  National Institute of Science Education and Research Bhubaneswar,  Bhubaneswar, Padanpur, Odisha 752050, India}

	
	\date{\today}
	\keywords{ Hypergraphs; Hypergraphs and linear algebra; Symmetries in hypergraphs; Invariant subspaces of hypergraph matrices; Eigenvalues and eigenvectors of hypergraph matrices. 
	}
	
	\subjclass[2020]{Primary
		05C65 
  05C50, 
  05C69, 
  05C70 
  	; Secondary
  05C75,
		05C07,
		05C30
	}
	
	\doublespacing
	 \maketitle
 \hrule
 \vspace{10pt}
	\begin{abstract}
 Here, the structural symmetries of a hypergraph are represented through equivalence relations on the vertex set of the hypergraph.
  A matrix associated with the hypergraph may not reflect a specific structural symmetry. In the context of a given symmetry within a hypergraph, we investigate a collection of matrices that encapsulate information about the symmetry.  Our investigation reveals that certain structural symmetries in a hypergraph manifest observable effects on the eigenvalues and eigenvectors of designated matrices associated with the hypergraph. 
  We identify specific matrices where the invariance is a consequence of symmetries present in the hypergraph. These invariant subspaces elucidate analogous behaviours observed in certain clusters of vertices during random walks and other dynamical processes on the hypergraph.

	\end{abstract}
 
 \vspace{10pt}
	\hrule
 
	\section{Introduction} 	
 The relationship between the structural symmetries of graphs and the spectra of matrices associated with graphs has been well-studied. These symmetries are represented by graph homomorphisms, graph coverings, quotient graphs, twin vertices, graph automorphisms, and quotient matrices. Effects of these symmetries in the spectra of matrices related to the graph are explored in \cite{Steve-Butler-2010-edgecover,Trevisan_symm_simax_2016,Oscar-rojo-balanced-tree-2005,mehatari-banerjee-motif2015,Vladimir-R-Covering-automorphisms,Benjamin-webb-gen-equitable-2019}.
  Often, these symmetries lead to invariant subspaces of matrices associated with graphs. For example, given any automorphism $\rho$ of a graph $G$, the matrix representation $P$ of $\rho$ commutes with the adjacency matrix $A_G$ of the graph $G$, that is $A_GP=PA_G$ (see \cite[Proposition 3.8.1]{rowlishionn-simic-sgt}). Consequently, given any eigenvalue $\lambda$ of $P$, the eigenspace of $\lambda$ is an invariant subspace of $A_G$. Equitable partitions of a graph, which is a generalisation of the notion of a graph symmetry,  also lead to invariant subspaces. Recall the following result from  \cite{agt-godsil-royle}.
  \begin{thm}(\cite[Lemma 9.3.1]{agt-godsil-royle})  Let $G$ be a graph with adjacency matrix $A_G$.
      If $\pi$ is an equitable partition of the graph $G$ with characteristic
matrix $S$, and $B$ is the adjacency matrix $A_{G/\pi}$ of the divisor
graph $G/\pi$, then $A_GS=SB$, and $B=(S^TS)
^{-1}S^TA_GS$.
  \end{thm}
  The equality $A_GS=SB$ means the column space $Col(S)$ of the characteristic
matrix $S$ is an invariant subspace of the adjacency matrix $A_G$ (see \cite[Lemma 9.3.2 ]{agt-godsil-royle}). This invariant subspace has notable implications. If $V(G)$ is the vertex set of $G$ and $V(G)=V_1\sqcup V_2\sqcup\ldots\sqcup V_m$ is the equitable partition $\pi$, then the $k$-th column of the characteristic matrix $S$ is the characteristic vector $\chi_{V_k}$ of $V_k$. That is, for some $i\in V(G)$ the $k$-th column $\chi_{V_k}$ has $i$-th entry $1$ if $i\in V_k$, otherwise the $i$-th entry is $0$. For any vector $x=\sum\limits_{i=1}^mc_i\chi_{V_i}\in Col(S)$, if $i,j$ are any pair of vertices
belonging to the same cell $V_k$, the $i$-th, and $j$-th components of $x$ are equal. Since $Col(S)$ is an invariant subspace of $A_G$, the equality of $i$-th and $j$-th components is preserved in $A_Gx(\in Col(S))$. Thus, if a vector $x$ is constant on each  $V_k$, this property is preserved under the action of the graph adjacency, and $A_Gx$ is also constant in $V_k$. That is, the symmetry encoded by the equitable partition induces an equality-preserving effect through the adjacency matrix. 

Hypergraphs are a generalisation of graphs in which the notion of edges is replaced by the notion of hyperedges. Instead of being a two-element subset of the vertex set, a hyperedge can be any non-empty subset of the vertex set. Like graphs, hypergraphs are also studied using various matrices, such as adjacency matrix, Laplacian matrix, signless Laplacian matrix etc (\cite{Banerjee-2021-hgmat,Sarkar-banerjee-2020,up2021,rodriguez2003laplacian,Swarup-panda-2022-hypergraph,cardoso-trevisan-2022-signless}). Consequently, a pertinent question arises: similar to graphs, do the symmetries of hypergraphs yield invariant subspaces of the matrices associated with them? We explore the answer to this question in this article.

Formally, A \emph{hypergraph} $H$ is an ordered pair $(V(H), E(H))$ such that $V(H)$ is a non-empty set and $E(H)$ is a family of non-empty subsets of $V(H)$. Each element of $V(H)$ is called a \emph{vertex} in $H$, and $V(H)$ is called the \emph{vertex set} of $H$. Any element of $E(H)$ is called a \emph{hyperedge}, and $E(H)$ is called the \emph{hyperedge set} in $H$.
For any vertex $i\in V(H)$, the \emph{star of the vertex} $i$ is $E_i(H) = \{e \in E(H) : i \in e\}$. For graphs, there are several ways to represent structural symmetries, such as graph automorphisms, equitable partitions, etc. A natural question is how the notion
of structural symmetry can be extended to hypergraphs.
For example, in \Cref{fig:unit}, the pair of vertices in $6,8\in V(H)$ have the same stars. Thus, $U=\{6,8\}$ is a symmetric substructure in $H$, as these vertices form an
equivalence class on $V(H)$. As this suggests, it is possible to represent this type of symmetry by using equivalence relations on the set of vertices of the hypergraph (see \cite{unit}). Given any equivalence relation $\mathfrak{R}$ on the vertex set, all the vertices belonging to an $\mathfrak{R}$-equivalence class are similar in some sense. Thus, $\mathfrak{R}$ represents a symmetry of the hypergraph.  Since any matrix associated with a hypergraph $H$ may contain only partial information about $H$, given a particular equivalence relation $\mathfrak{R}$ on the vertex set $V(H)$, a matrix related to $H$ may or may not encode the information about the symmetry represented by $\mathfrak R$. Thus, it is natural to ask what family of matrices associated with a hypergraph contains the information about the equivalence relation. In \cite{up2021}, we introduced some general operators related to hypergraphs. 
 {In \cite{unit}, for any equivalence relation on $V(H)$, we have introduced the notion of equivalence relation compatible operators, a family of operators that can encrypt specific information of the symmetry represented by the equivalence relation. In that paper we have studied some symmetries of hypergraphs, which are captured by the general operators associated with the hypergraph. In this article, we find some other classes of matrices that can encrypt information about hypergraph symmetries in terms of their invariant subspaces and as well as their eigenvalues. 
  These families of matrices encompass certain adjacency and Laplacian matrices associated with hypergraphs that have previously been studied in the literature (see \Cref{ex:adj} and \Cref{ex:lap}). However, they also lead to various new matrices related to hypergraphs. For example, a \emph{star-dependent matrix} \( M_H \) (see \Cref{star-matrix}), associated with a hypergraph \( H \), is a square matrix whose rows and columns are indexed by the vertex set \( V(H) \). This matrix is defined based on two functions, \( g:\mathfrak{E}(H) \times \mathfrak{E}(H) \to \mathbb{R} \) and \( h:\mathfrak{E}(H) \to \mathbb{R} \), where \( \mathfrak{E}(H) = \{E_i(H) : i \in V(H)\} \) denotes the collection of all \emph{stars} in \( H \). Specifically, for any \( i \in V(H) \), the diagonal entry \( (i,i) \) of \( M_H \) is \( h(E_i(H)) \), while for distinct \( i, j \in V(H) \), the off-diagonal entry \( (i,j) \) is \( g(E_i(H), E_j(H)) \). Each choice of the pair \( (g, h) \) defines a different star-dependent matrix.
Certain known adjacency, Laplacian, and signless Laplacian matrices are star-dependent matrices arising from specific choices of \( (g, h) \). However, many other choices yield new star-dependent matrices. For instance, to the best of our knowledge, the star-dependent matrix presented in \Cref{finer-ex}(1) appears to be a new hypergraph matrix in the literature.
Among the star-dependent matrices, \emph{edge-sum matrices} form a notable subclass where \( (g, h) \) depends on edge functions \( p_1: E(H) \to \mathbb{R} \) and \( p_2: E(H) \to \mathbb{R} \). In these cases, particular choices of \( (p_1, p_2) \) correspond to well-known hypergraph matrices, while many other choices introduce new types of hypergraph matrices.

 Certain symmetries inherent in hypergraphs give rise to invariant subspaces of matrices associated with the hypergraph.
   For example, here we find star-dependent matrices that encode the information about the symmetries represented by the equivalence relation $\mathfrak{R}_s(H)$ that depends on the star of the vertices. For two vertices $i$, and $j$, the pair  $(i,j)\in\mathfrak{R}_s(H)$, if they have the same star. The equivalence relation  $\mathfrak{R}_2(H,M_H)$ is related to the hypergraph $H$ and a matrix $M_H$ associated with $H$. A pair of vertices $(i,j)$ is called $\mathfrak{R}_2(H,M_H)$ related if there exists a bijection between the stars of $i$ and $j$, and the symmetry represented by the bijection is encrypted in the matrix $M_H$. We show that if $M_H$ belongs to a specific class of matrices, referred to as edge-sum matrices, then $\mathfrak{R}_2(H,M_H)$-equivalence class leads to invariant subspaces of $M_H$.  These symmetries, represented by the equivalence relations  $\mathfrak{R}_s(H)$ and $\mathfrak{R}_2(H,M_H)$ produce some eigenvalues and eigenvectors of the star-dependent matrices and edge-sum matrices, respectively. 
 We show that these eigenvectors provide invariant subspaces of related matrices.  
 These invariant subspaces preserve the equality of dynamics in some clusters of vertices during dynamical processes in hypergraphs. That is, once the vector representing the state of the dynamics becomes equal in a cluster of vertices, the equality is preserved in subsequent times.

  \Cref{sec-inv} starts with an equivalence relation $\mathfrak{R}_{M_H}$ on the vertex set $V(H)$ of a hypergraph $H$. This equivalence relation represents the symmetries inherent in the matrix $M_H$.

 We show that certain eigenvalues and eigenvectors of \( M_H \) can be expressed in terms of the equivalence classes of \( \mathfrak{R}_{M_H} \) (see \Cref{eigen-lemma-1}). The eigenvectors corresponding to an \( \mathfrak{R}_{M_H} \)-equivalence class \( W \) form a basis for an invariant subspace \( S_W \) of \( M_H \). Furthermore, if \( M_H \) is symmetric with respect to a given inner product, then the orthogonal complement of \( S_W \) with respect to this inner product is also an invariant subspace of \( M_H \) (see \Cref{class_const}).
 The equivalence relation $\mathfrak{R}_{M_H}$ is related to the matrix $M_H$.  
Using  \Cref{finer_class_const} and  $\mathfrak{R}_{M_H}$,  we find that the structural symmetries of a hypergraph encoded by the equivalence relations $\mathfrak{R}_s(H)$, and $\mathfrak{R}_2(H,M_H)$ causes specific eigenvalues, eigenvectors of some family of matrices namely \emph{star-dependent matrices, edge-sum matrices, edge-cardinality matrices}. These eigenvalues and eigenvectors sometimes lead to invariant subspaces of these matrices. 
Each $\mathfrak{R}_s(H)$-equivalence class is called a unit in the hypergraph $H$. 
  A {star dependent} matrix (defined in \Cref{star-matrix}) is a matrix of order $|V(H)|$, whose $i,j$-th element depends only on the star of $i$, and $j$. We show that if $M_H$ is a star-dependent matrix, then the structural symmetries in the units lead us to invariant subspaces of $M_H$ (see \Cref{unit-finer}, \Cref{unit-invariant}, \Cref{unit_invariant_cluster}).
  An edge-sum matrix is a star-dependent matrix of order $|V(H)|$ such that for all $i,j\in V(H)$, the $i,j$-th entry of the matrix is a weighted sum on the intersection of stars, $E_i(H)\cap E_j(H)$. The structural symmetry of $H$ represented by each $\mathfrak{R}_2(H,M_H)$-equivalence class is manifested as invariant subspaces of edge-sum matrices associated with $H$ (see \Cref{twin-finer}, and \Cref{edge-sum-cor}). During random walks and other dynamical processes on hypergraphs, sometimes the dynamics on some clusters of vertices exhibit a tendency to become very similar to each other. We use our results to explain this behaviour in \Cref{sec-application}.
\section{Invariant Subspaces Due to Transpositions}\label{sec-inv}
 Let $H$ be a hypergraph and  $M_H$ be a real square matrix whose rows and columns are {indexed by} the set of vertices $V(H)$ in $H$, that is $M_H=\left(m_{ij}\right)_{i,j\in V(H)}$ with $m_{ij}\in\mathbb{R}$ for all $i,j\in V(H)$. 
		For any square matrix $M$, 
  we denote the set of all eigenvalues of $M$ as $\sigma_M$. Let $V_\lambda(M)$ be the eigenspace of $\lambda\in\sigma_M$, that is, $V_\lambda(M)$ is the vector space generated by $\{x:Mx=\lambda x\}$.
	
	For two distinct $i,j\in V(H)$, the transposition $\phi_{ij}:V(H)\to V(H)$ is defined as 
	\[\phi_{ij}(k)=
	\begin{cases}
		j &\text{~if~} k=i,\\
		i&\text{~if~}k=j,\\
		k&\text{~otherwise.}
	\end{cases}
	\]
	The bijection $\phi_{ij}$ leads us to the map $\psi_{ij}:\mathbb{R}^{V(H)}\to \mathbb{R}^{V(H)}$ defined by $\psi_{ij}(x)=x\circ\phi_{ij}$. 
 A subspace $W$ of the vector space $\mathbb{R}^{V(H)}$  is called an \emph{invariant} subspace of $M_H$ if $M_Hy\in W$ for all $y\in W$.
{ If two square matrices, $A$ and $B$ commute, the eigenspace of any eigenvalue of $A$ is invariant under the action of $B$.  As an example, any eigenspace of the adjacency matrix of a graph $G$ is invariant under the action of any automorphism of $G$  ( \cite[Section-1]{chan1997symmetry}). This observation motivates us to consider the following relation on $V(H)$.}	
\[\mf{R}_{M_H}:=\{(i,j)\in V(H)\times V(H):m_{ii}=m_{jj},m_{ij}=m_{ji}, m_{ik}=m_{jk},m_{ki}=m_{kj}\text{~for all~}k\notin\{i,j\}\}.\]
We say that vertices \( i \) and \( j \) are \(\mathfrak{R}_{M_H}\)-related if \((i, j) \in \mathfrak{R}_{M_H}\). This implies that the \(2 \times 2\) submatrix of \(M_H\) corresponding to the vertices \(i\) and \(j\) satisfies the following conditions: the diagonal entries are equal, and the off-diagonal entries are also equal. Furthermore, outside this \(2 \times 2\) submatrix, the rows and columns of \(M_H\) corresponding to \(i\) and \(j\) are identical.
From the definition of $\mf{R}_{M_H}$, it follows that $\mf{R}_{M_H}$ is reflexive and symmetric. If for three distinct $i,j,k\in V(H)$, we have $(i,j),(j,k)\in \mf{R}_{M_H}$, then  $m_{ii}=m_{jj}=m_{kk}$, $m_{ik}=m_{jk}=m_{kj}=m_{ki}$, and for all $l\notin\{i,j,k\}$ we have $m_{il}=m_{jl}=m_{kl} $ and $m_{ij}=m_{ji}=m_{ki}=m_{kj}$. That is, $m_{il}=m_{kl}$ for all $l\notin\{i,k\}$. Similarly, we can show $m_{li}=m_{lk}$ for all $l\notin\{i,k\}$. Therefore, $(i,j),(j,k)\in \mf{R}_{M_H}$ implies $(i,k)\in \mf{R}_{M_H}$, that is $\mf{R}_{M_H}$ is transitive. Thus, $\mf{R}_{M_H}$ is an equivalence relation. Suppose that $W\subseteq V(H)$ is an $M_H$-equivalence class. If four distinct vertices $i,j,k,l\in W$, then $ m_{ij}=m_{ji}=m_{ki}=m_{kl}$. That is, each $\mathfrak{R}_{M_H}$-equivalence class $W$ corresponds to two constants $d_W$ and $c_W$ such that $W$ is a maximal subset of $V(H)$ with $m_{ii}=d_W$, and $m_{ij}=c_W$ for all $i,j(\ne i)\in W$. Consequently, $W$ corresponds to a $|W|\times |W|$ submatrix of $M_H$ of the form $d_WI_{|W|}+c_W(J_{|W|}-I_{|W|})$, where $I_n$ is the identity matrix of order $n$, and $J_n$ is an $n\times n$ matrix with all entries equal to 1.

 Now, we show  given any $(i,j)\in \mathfrak{R}_{M_H}$, the matrix $M_H$ commutes with the $\psi_{ij}$. Like the commutativity of the adjacency matrix of a graph
 and the matrix representation of an automorphism of the graph, this commutativity also leads us to $M_H$-invariant subspaces. For $i,j,k,k'\in V(H)$, with $(i,j)\in \mathfrak{R}_{M_H}$, if $k\notin\{i,j\}$, then $m_{ki}x(i)+m_{kj}x(j)=m_{kj}(x\circ \phi_{ij})(j)+m_{ki}(x\circ \phi_{ij})(i)$. Consequently,
  \begin{align*}
    & ((M_Hx)\circ\phi_{ij})(k)=  (M_Hx)(k)=\sum_{k'\in V(H)}m_{kk'}x(k')=\sum_{k'\in V(H)}m_{kk'}(x\circ \phi_{ij})(k')=(M_H(x\circ\phi_{ij}))(k).
     \end{align*}
  Now, $m_{ij}=m_{ji}$, and $m_{ii}=m_{jj}$, lead us to $ m_{ii}x(i)+m_{ij}x(j)=m_{jj}x(i)+m_{ji}x(j)=m_{jj}(\phi_{ij}\circ x)(j)+m_{ji}(\phi_{ij}\circ x)(i)$.   Therefore, $m_{ik}=m_{jk}$ for all $k\notin \{i,j\}$ implies that
    \begin{align*}
        &((M_Hx)\circ\phi_{ij})(j)=(M_Hx)(i)=\sum_{k\in V(H)}m_{ik}x(k)=\sum_{k\in V(H)}m_{jk}(x\circ\phi_{ij})(k)=(M_H(x\circ\phi_{ij}))(j).
    \end{align*}
 Proceeding similarly, we can show that $((M_Hx)\circ\phi_{ij})(i) =(M_H(x\circ\phi_{ij}))(i)$. Thus, for all $(i,j)\in \mathfrak{R}_{M_H}$,
 $$M_H(\psi_{ij}(x))=\psi_{ij}(M_H x)\text{~for~all}~x\in {\mathbb{R}}^{V(H)}.$$

\begin{exm}\label{exm-MH}
   Consider the hypergraph \( H \) with vertices \( V(H) = \{1, 2, 3, 4\} \) and edges \( E(H) = \{e, f, g\} \), where \( e = \{1, 2, 4\} \), \( f = \{2, 3, 4\} \), and \( g = \{1, 3, 4\} \). Define the matrix \( M_H = (m_{ij})_{i,j \in V(H)} \) associated with \( H \), by \( m_{ij} = \sum\limits_{e \in E_i(H) \cap E_j(H)} \frac{1}{|e|} \) for all \( i, j \in V(H) \). The matrix \( M_H \) forms two \( \mathfrak{R}_{M_H} \)-equivalence classes: \( \{1, 2, 3\} \) and \( \{4\} \).
\begin{figure}[ht]
  \begin{center}
     \begin{tikzpicture}[scale=0.5]
		\node [style=new style 0,scale=0.6] (0) at (0, 0) {$4$};
		\node [style=new style 0,scale=0.6] (1) at (-2, -2) {$2$};
		\node [style=new style 0,scale=0.6] (2) at (2, -2) {$3$};
		\node [style=new style 0,scale=0.6] (3) at (0, 2.5) {$1$};
		\node [style=none] (4) at (-2.75, -2.75) {};
		\node [style=none] (5) at (0, 3.5) {};
		\node [style=none] (6) at (1.25, -0.5) {};
		\node [style=none] (7) at (2.25, -2.75) {};
		\node [style=none] (8) at (-0.5, 3.5) {};
		\node [style=none] (9) at (-1.75, -0.5) {};
		\node [style=none] (10) at (0, 1.5) {};
		\node [style=none] (11) at (-3, -2.25) {};
		\node [style=none] (12) at (2.5, -2.25) {};
		\node [style=none] (13) at (-3, 1.5) {e};
		\node [style=none] (14) at (-0.25, -4.25) {f};
		\node [style=none] (15) at (2.5, 1.5) {g};
	    \node [style=none] (16) at (-2, -4.25) {H};
            \node [style=none] (17) at (9, 0.5) {$M_H=$
           
\begin{tabular}{ c||ccc|c| }

               &1&2&3&4\\
               \hline\hline
              1&  $\frac{2}{3}$&$\frac{1}{3}$&$\frac{1}{3}$&$\frac{2}{3}$\vspace{3pt}\\
              2&  $\frac{1}{3}$&$\frac{2}{3}$&$\frac{1}{3}$&$\frac{2}{3}$\vspace{3pt}\\
               3&$\frac{1}{3}$&$\frac{1}{3}$&$\frac{2}{3}$&$\frac{2}{3}$\vspace{3pt}\\
               \hline
              4&  $\frac{2}{3}$&$\frac{2}{3}$&$\frac{2}{3}$&$1$\vspace{3pt}\\  \hline
\end{tabular}};
		\draw (4.center)
			 to [bend right=60, looseness=0.50] (6.center)
			 to [bend right, looseness=0.50] (5.center)
			 to [bend left=285, looseness=0.75] cycle;
		\draw (7.center)
			 to [bend left=60, looseness=0.50] (9.center)
			 to [bend left, looseness=0.50] (8.center)
			 to [bend left=75, looseness=0.75] cycle;
		\draw (12.center)
			 to [bend right=45, looseness=0.50] (10.center)
			 to [bend right=45, looseness=0.50] (11.center)
			 to [bend right=45, looseness=1.25] cycle;
	
\end{tikzpicture}
    \caption{A matrix $M_H$ associated with hypergraph $H$ and the relation $\mathfrak{R}_{M_H}$.}
    \label{fig:MH}
  \end{center}
\end{figure} 
   Here $(1,2)\in\mathfrak{R}_{M_H}$. For any $x\in\mathbb{R}^{V(H)}$,  we have $(M_Hx)(1)=\frac{1}{3}(2x(1)+x(2)+x(3)+2x(4))$, and  $(M_Hx)(2)=\frac{1}{3}(x(1)+2x(2)+x(3)+2x(4))$. Thus in $x$, if we interchange $x(1)$ with $x(2)$, then the values of $(M_Hx)(1)$ and $(M_Hx)(2)$ get interchanged. That is, $M_H(\psi_{12}(x))(1)=(M_H(x))(2)=\psi_{12}(M_H x)(1)\text{~for~all}~x\in {\mathbb{R}}^{V(H)}.$

\end{exm}
 
Any positive real-valued function $f:V(H)\to (0,\infty)$ provides an inner product $(\cdot,\cdot)_f$ on $\mathbb{R}^{V(H)}$,  defined as $$(x,y)_f=\sum\limits_{i\in V(H)}f(i)x(i)y(i)\text{~for all~} x,y\in \mathbb{R}^{V(H)}.$$
	 We say $M_H$ is an $f$-symmetric matrix\footnote{In most of the results, we assume $M_H$ to be $f$-symmetric. This enables us to consider the eigenvectors of $M_H$ as real eigenvectors in $\mathbb{R}^{V(H)}$. However, in \Cref{unit_invariant_cluster}, $M_H$ is not assumed to be $f$-symmetric, and the eigenvectors of $M_H$ are not used to prove the theorem.
   } if $(M_Hx,y)_f=(x,M_Hy)_f$ for all $x,y\in\mathbb{R}^{V(H)}$.
	For any $x\in \R^{V(H)}$, the support of $x$ is $supp(x)=\{i\in V(H):x(i)\ne 0\}$. We define the following for any $U\subseteq V(H)$.
 
 (1) $S_U=\{x\in \mathbb{R}^{V(H)}: supp(x)\subseteq U, \sum\limits_{i\in U}x(i)=0\}$.

 (2) For any $f:V(H)\to (0,\infty)$, $C^f_U=\{x\in\mathbb{R}^{V(H)}:f(i)x(i)=c_x,\text{~a real constant for all~}i\in U\}$. If $f(i)=1$ for all $i\in V(H)$, then we denote $C^f_U$ as $ C_U=\{x\in \mathbb{R}^{V(H)}:x(i)=c_x\text{~for all~}i\in U\}$.

 (3) The characteristic function $\chi_{U}:V(H)\to\mathbb{R}$ of $U$ is defined as $\chi_{_U}(i)=\begin{cases}
			1&\text{~if~}i\in U,\\
			0&\text{otherwise.}
		\end{cases} $
	For two distinct $i,j\in V(H)$, the function $x_{ij}=\chi_{\{j\}}-\chi_{\{i\}}$ is referred to as a \emph{Faria vector} (see \cite[p.46, Section 3.7]{biyikoglu2007laplacian}).
	For any $U=\{i_0,i_1,\ldots,i_{n-1}\}\subseteq V(H)$, the collection $\{x_{i_ji_0}:j=1,\ldots,n-1\}$ forms a basis of $S_U$. Therefore, $S_U$ is a subspace of dimension $|U|-1$. The next result shows a relation between $S_U$ and $C^f_U$. For any subspace $\mathcal{V}$ of $\mathbb{R}^{V(H)}
	$, we denote the \emph{orthogonal complement\footnote{A matrix \( M_H \) is \( f \)-symmetric if \((M_Hx, y)_f = (x, M_Hy)_f\) for all vectors \( x \) and \( y \) in \(\mathbb{R}^{V(H)}\). This means that if \(\mathcal{V}\) is a subspace invariant under \( M_H \), then its orthogonal complement \(\mathcal{V}^{(\perp,f)}\) is also invariant under \( M_H \). This fact is used to prove \Cref{class_const}} of $\mathcal{V}$ with respect to the inner product $(\cdot,\cdot)_f $} to be $ \mathcal{V}^{(\perp,f)}$.
	\begin{prop}\label{perp_lem}
		Let $H$ be a hypergraph. For any $U\subseteq V(H)$ with $ |U|>1$, and any function $f:V(H)\to (0,\infty)$,
		the orthogonal complement of $S_U$ with respect to the inner product $(\cdot,\cdot)_f$, $S_U^{(\perp,f)}=C^f_U$.
	\end{prop}
	\begin{proof}Let $U=\{i_0,i_1,\ldots,i_{n-1}\}\subseteq V(H)$.
		For any $y\in S_U^{(\perp,f)}$, if $i,j(\ne i)\in U$, then $$\sum_{u\in V(H)}f(u)y(u)x_{ij}(u)=0,$$ and so $f(i)y(i)=f(j)y(j)$ for all $i,j(\ne i)\in U$. Therefore, $y\in C^f_U$, and thus, $S_U^{(\perp,f)}\subseteq C^f_U$. Conversely, for any $y\in C_U^f$, $f(i)y(i)=c_y$ for all $i\in U$. Thus, if $x\in S_U$, then $$\sum\limits_{i\in V(H)}f(i)x(i)y(i)=\sum\limits_{i\in U}f(i)x(i)y(i)=c_y\sum\limits_{i\in U}x(i)=0,$$ and therefore, $y\in S_U^{(\perp,f)}$. Thus, $C^f_U\subseteq S_U^{(\perp,f)}$ and the result follows.
	\end{proof}
	For $(i,j)\in \mathfrak{R}_{M_H}$ a symmetry between $i$ and $j$ in $H$ is encoded in $M_H$ in such a way that the action of $M_H$ commutes with the transposition between the vertices $i$ and $j$. In the next result, we show that this symmetry contributes to the spectra of $M_H$.
	\begin{prop}\label{eigen-lemma-1}
		Let $H $ be a hypergraph. For any $\mf{R}_{M_H}$-equivalence class $W$ with $|W|\ge 2$, there exists $\lambda_W\in \sigma_{M_H}$, such that $S_W$ is a subspace of $V_{\lambda_W}(M_H)$.
	\end{prop}
	\begin{proof}
		Suppose that $W=\{i_0,i_1,\ldots,i_{n-1}\}$. For any $i,j\in W$ we have $\phi_{ij}(x_{ij})=-x_{ij}$. Thus, if $u\in V(H)\setminus\{i,j\}$ then 
		\begin{align*}
			& (\phi_{ij}(M_H x_{ij}))(u)=M_H(\phi_{ij}( x_{ij}))(u)\\
			&\implies (M_H x_{ij})(u)=0,
		\end{align*}
		and $(\phi_{ij}(M_H x_{ij}))(i)=(M_H(\phi_{ij}( x_{ij})))(i)$ implies $(M_H x_{ij})(j)=-(M_H x_{ij})(i)$. Therefore, $M_Hx_{ij}=-[(M_Hx_{ij})(i)]x_{ij}$.
  
  If \( |W| = 2 \), then the dimension of \( S_W \) is \( |W| - 1 = 1 \). Since \( x_{i_0 i_1} \) is a non-zero vector in \( S_W \), it follows that the one-dimensional vector space \( S_W \) is spanned by \( \{ x_{i_0 i_1} \} \). Consequently, \( \lambda_W = -[(M_H x_{i_0 i_1})(i_0)] \) is an eigenvalue of \( M_H \), and \( S_W \) forms a subspace of \( V_{\lambda_W}(M_H) \).
  
  If $|W|>2$, then
		for three distinct $i,j,k\in W $, since $(M_H(\phi_{jk}(x_{ij})))(i)=(\phi_{jk}( M_H(x_{ij})))(i)$, we have $ (M_Hx_{ik})(i)=( M_Hx_{ij})(i)$. Therefore, $-[(M_Hx_{i_0i_j})(i_0)]=\lambda_W$, a constant for all $j=1,\ldots,n-1$ and, 
		
		\begin{equation}\label{eigeneqn}
			M_Hx_{i_0i_j}=\lambda_Wx_{i_0i_j}.
		\end{equation}    
		Since $\{x_{i_0i_j}:j=1,\ldots,n-1\}$ is a basis for $S_W$, we have $S_W$ is a subspace $V_{\lambda_W}(M_H)$. This completes the proof.
	\end{proof}
We demonstrate \Cref{eigen-lemma-1} in the following example.
\begin{exm}
    Recall the hypergraph $H$ and the matrix $M_H$ are considered in  \Cref{exm-MH} (see \Cref{fig:MH}). For $W=\{1,2,3\}$,
    consider the subspace $S_W=\{x\in\mathbb{R}^{V(H)}:x(4)=0,\text{~and~}x(1)+x(2)+x(3)=0\}$.
    The subspace $S_W$ can be written as $S_W=\{(x(1),x(2),-(x(1)+x(2)),0):x(1),x(2)\in \mathbb{R}\}=\{x(1)(1,0,-1,0)^T+x(2)(0,1,-1,0)^T:x(1),x(2)\in \mathbb{R}\}=\{x(1)x_{31}+x(2)x_{32}:x(1),x(2)\in \mathbb{R}\}$. Therefore, $S_W$ is a $2$-dimensional vector space generated by $x_{31},x_{32}$.
    Again for any $i\in V(H)$, we have $ M_Hx_{31}=M_H(1,0,-1,0)^T=(m_{11}-m_{13},m_{21}-m_{23},m_{31}-m_{33},m_{41}-m_{43})^T$. Since $W=\{1,2,3\}$ is an $\mathfrak{R}_{M_H}$-equivalence class, we have $m_{21}-m_{23}=0=m_{41}-m_{43}$ and $m_{11}-m_{13}=-(m_{31}-m_{33}) $. Therefore,  $ M_Hx_{31}=(m_{11}-m_{13})x_{31}$. Similarly, we can show $ M_Hx_{32}=(m_{22}-m_{23})x_{32}$ and $m_{11}-m_{13}=m_{22}-m_{23}$. Therefore the vector space $S_W$, spanned by $x_{31}$ and $x_{32}$, is a subspace of the eigenspace $V_{m_{11}-m_{13}}(M_H)$. Now if $M_H$ is the same matrix that is considered in \Cref{exm-MH}, then $m_{11}-m_{13}=\frac{1}{3} $ and $S_W$ is a subspace of  $V_{m_{11}-m_{13}}(M_H)=V_{\frac{1}{3}}(M_H)$.
\end{exm}
	
	In the following Theorem, for any matrix $M_H$, each $\mf{R}_{M_H}$-equivalence class with cardinality at least $2$ corresponds to an invariant subspace of $M_H$.
	\begin{thm}\label{class_const}
		Let $H$ be a hypergraph, $f:V(H)\to (0,\infty)$ be a function, and $M_H$ be an $f$-symmetric matrix. For any $\mf{R}_{M_H}$-equivalence class $W$ with $|W|\ge 2$, the vector space $C_W^f$ is an invariant subspace of $M_H$.
	\end{thm}
	\begin{proof}
		By \Cref{eigen-lemma-1}, for an $\mf{R}_{M_H}$-equivalence class $W$ with $|W|\ge 2$, there exists an eigenvalue $\lambda_W$ of $M_H$, such that, $S_W$ is a subspace of the eigenspace $V_\lambda(M_H)$. Thus, $S_W$ is an invariant subspace of $M_H$. Since by \Cref{perp_lem} $C_W^f=S_W^{(\perp,f)}$, and $M_H$ is $f$-symmetric, the vector space $C_W^f$ is an invariant subspace of $ M_H$.
	\end{proof}
 A vector $x:V(H)\to\mathbb{R}$ is said to be constant on a subset $U\subseteq V(H)$ if $x(v)=c$, a fixed real for all $v\in U$.
	\Cref{class_const} indicates that for any hypergraph $H$ and any matrix $M_H$ indexed by $V(H)$, if $\{x_n:V(H)\to\mathbb{R}:n\in\mathbb{N}\}$ is governed by the difference equation
	\begin{equation}\label{mheqn}
		x_{n+1}=M_Hx_n,
	\end{equation}
	and there exists $k\in \mathbb{N}$ such that $x_k$ is constant in a $\mathfrak{R}_{M_H}$-equivalence class $W$ then so is $x_n$ for all $n\ge k$. Since many phenomena on hypergraphs can be described using difference equations similar to  \Cref{mheqn}, in the remaining part of this article, we use \Cref{class_const} to explore the effect of symmetries of hypergraphs on these phenomena. For example, consider the hypergraph $H$ and the matrix $M_H$ described in \Cref{exm-MH}. Suppose that $x_1:V(H)\to\mathbb{R}$ with $x_1(1)=x_1(2)=x_1(3)=c_1$, a real constant. Given $x_2=M_Hx_1$, we have $x_2(i)=(m_{i1}+m_{i2}+m_{i3})c_1+m_{i4}x_1(4)$. Since $W=\{1,2,3\}$ is an $\mathfrak{R}_{M_H}$ equivalence class, $m_{ij}=m_{i'j}$ for all $i,i'\in W$. Thus, the $i$th entry $x_2(i)=(m_{i1}+m_{i2}+m_{i3})c_1+m_{i4}x_1(4)=c_2$, another constant for all $i\in W$. Similarly, for any $n\ge 2$, the vector $x_n$ would be a constant on $W$. Thus, the symmetry of $W$ in $H $ is reflected in the fact that $x_n(1)=x_n(2)=x_n(3)$ for all $n\ge 2$. We will revisit \Cref{mheqn} once we have proven \Cref{unit_invariant_cluster}.

As described in the definition of $\mathfrak{R}_{M_H}$ and the discussion that follows it, the equivalence relation $\mathfrak{R}_{M_H}$ is related to the matrix $M_H$. Now, we relate some structural symmetries of hypergraphs with $\mathfrak{R}_{M_H}$.  
Let $\mathfrak{R}_1$ and $\mathfrak{R}_2$ denote two equivalence relations on the vertex set $V(H)$. The relation $\mathfrak{R}_1$ is considered \emph{finer than} $\mathfrak{R}_2$ if $ \mathfrak{R}_1\subseteq \mathfrak{R}_2$. This implies that if $(i,j)\in \mathfrak{R}_1$, then $(i,j)\in \mathfrak{R}_2$, indicating that for any $\mathfrak{R}_1$-equivalence class $W_1$, there exists a corresponding $\mathfrak{R}_2$-equivalence class $W_2$ such that $W_1\subseteq W_2$.

\begin{exm}\label{ex-finer-1}
    Let $H$ be a hypergraph. 
    Consider the relation $\mathfrak{R}_s(H)=\{(i,j)\in V(H)\times V(H):E_i(H)=E_j(H)\}$ on the vertex set $V(H)$. Since $(i,i)\in \mathfrak{R}_s(H)$ for all $i\in V(H)$, the equivalence relation $\mathfrak{R}_s(H)$ is reflexive. It is symmetric because $(i,j)\in \mathfrak{R}_s(H)$ means $E_i(H)=E_j(H)$ which implies $(j,i)\in \mathfrak{R}_s(H)$. If $(i,j),(j,k)\in  \mathfrak{R}_s(H)$ then $E_i(H)=E_j(H)=E_k(H) $, that is $(i,k)\in \mathfrak{R}_s(H)$. Consequently, $\mathfrak{R}_s(H)$ is an equivalence relation.   
    Consider the matrix \( M_H = (m_{ij})_{i,j \in V(H)} \) associated with \( H \), where \( m_{ij} = \sum\limits_{e \in E_i(H) \cap E_j(H)} \frac{1}{|e|} \) for all \( i, j \in V(H) \). For all $(i,j)\in \mathfrak{R}_s(H)$ we have $E_i(H)=E_j(H)$. Therefore, $m_{ii}=\sum\limits_{e\in E_i(H)}\frac{1}{|e|}=\sum\limits_{e\in E_j(H)}\frac{1}{|e|}=m_{jj}$. From the definition of $M_H$, it follows $m_{ij}=m_{ji}$. For any $k\notin\{i,j\}$, we have $m_{ik}=m_{ki}=\sum\limits_{e\in E
_i(H)\cap E_k(H)}\frac{1}{|e|}=\sum\limits_{e\in E
_j(H)\cap E_k(H)}\frac{1}{|e|}=m_{jk}=m_{kj}$. Therefore $(i,j)\in \mathfrak{R}_{M_H}$. That is, for this matrix $M_H$, the equivalence relation $\mathfrak{R}_s(H)$ is finer than the equivalence relation $\mathfrak{R}_{M_H}$. In the hypergraph $H$ is considered in \Cref{exm-MH}, all the $\mathfrak{R}_s(H)$ equivalence class are singleton sets, and the $\mathfrak{R}_{M_H}$ equivalence classes are $\{1,2,3\}$, and $\{4\}$. Thus, for all $\mathfrak{R}_s(H)$ equivalence class $W_1$, there exists a  $\mathfrak{R}_{M_H}$ equivalence class $W_2$ such that $W_1\subseteq W_2$.
Later in \Cref{unit-finer}, we provide a family of matrices such that $\mathfrak{R}_s$ is finer than the equivalence relation $\mathfrak{R}_{M_H}$ when $M_H$ belongs to that family.
\end{exm}

As demonstrated in \Cref{class_const}, it is established that $C^f_W$ is an invariant subspace of $M_H$, where $W$ represents an $\mathfrak{R}_{M_H}$-equivalence class. We now extend this result to assert that the same holds for any $\mathfrak{R}$-equivalence class $W$, where $\mathfrak{R}$ is an equivalence relation finer than $\mathfrak{R}_{M_H}$.
	
	\begin{thm}\label{finer_class_const}
		Let $H$ be a hypergraph, $f:V(H)\to (0,\infty)$ be a function, and $M_H$ be an $f$-symmetric matrix. If $\mathfrak{R}$, an equivalence relation on $V(H)$, is finer than $\mf{R}_{M_H}$ then for any $\mf{R}$-equivalence class $W$ with cardinality $|W|\ge 2$, the vector space $C_W^f$ is an invariant subspace of $M_H$.
	\end{thm}
	\begin{proof}
		Since for any  $\mf{R}$-equivalence class $W$, there exists an  $\mf{R}_{M_H}$-equivalence class $W'$ such that $W\subseteq W'$, and thus, $S_W$ is a subspace of $S_{W'}$. So, from \Cref{eigeneqn} we
		conclude that $S_W$ is a subspace generated by some eigenvectors of the eigenvalue $\lambda_{W'}\in \sigma_{M_H}$. Thus, $S_W$ is an invariant subspace of $M_H$, and so is $C_W^f$.
	\end{proof}

 Now, we introduce some equivalence relations related to structural symmetries of $H$, which are finer than $\mathfrak{R}_{M_H}$. We use \Cref{finer_class_const} to explore the relation of these structural symmetries of the hypergraph $H$ with the matrix $M_H$.
 
	\subsection{Units, Star Dependent Matrices and their Invariant Subspaces }
 
	\subsection*{Units in a Hypergraph} In a hypergraph $H$, \emph{units} are the maximal subsets of vertices having the same star.
	\begin{df}[Units, \cite{unit}]\label{unit}
		Let $H$ be the hypergraph and $\mathfrak{R}_s(H)$ be an equivalence relation on $V(H)$ defined by
		$$\mathfrak{R}_s(H)=\{(i,j)\in V(H)\times V(H):E_i(H)=E_j(H)\}.$$
		Then, each $\mathfrak{R}_s(H)$-equivalence  class in $H$ is called a \emph{unit}. 
	\end{df}
	A unit in $H$ corresponds to an $E\subseteq E(H)$, such that the unit is $\{i\in V(H):E_i(H)=E\}$, the maximum collection of vertices $i$, with $E_i(H)=E$. We denote the unit as $W_E$ and refer to $E$ as the \emph{generating set of the unit $W_E$.} We denote the collection of all the units in a hypergraph $H$ as $\mathfrak{U}(H)$.
	
	\begin{figure}[ht]
		\colorlet{cl1}{gray!20!white}
		\centering
		\begin{tikzpicture}
	
		\node [style=none] (0) at (-0.25, 0) {};
		\node [style=none] (1) at (1, 4) {};
		\node [style=none] (2) at (3.75, 1.25) {};
		\node [style=none] (3) at (1.25, 0) {};
		\node [style=none] (4) at (-0.25, 4) {};
		\node [style=none] (5) at (-2.5, 1.25) {};
		\node [style=none] (6) at (0, 1.25) {};
		\node [style=none] (7) at (2.25, -2.5) {};
		\node [style=none] (8) at (3.75, -0.25) {};
		\node [style=none] (9) at (1, 1) {};
		\node [style=none] (10) at (-0.75, -2.5) {};
		\node [style=none] (11) at (-2.75, -0.5) {};
		\node [style=new style 0,scale=0.7] (12) at (0.5, 0.5) {1};
		\node [style=new style 0,scale=0.7] (13) at (-2, -1) {2};
		\node [style=new style 0,scale=0.7] (14) at (-1, -1.75) {3};
		\node [style=new style 0,scale=0.7] (15) at (2, -1.5) {4};
		\node [style=new style 0,scale=0.7] (16) at (2.75, -0.5) {5};
		\node [style=new style 0,scale=0.7] (17) at (0.5, 2.25) {6};
		\node [style=new style 0,scale=0.7] (18) at (2.25, 2.5) {7};
		\node [style=new style 0,scale=0.7] (19) at (-1.25, 2.5) {9};
		\node [style=new style 0,scale=0.7] (20) at (-1.75, 1.5) {10};
		\node [style=none] (22) at (0.5, -2.75) {$H$};
		\node [style=none] (25) at (-0.25, 3.5) {};
		\node [style=none] (26) at (-2, 1) {};
		\node [style=none] (27) at (2.25, 3) {};
		\node [style=none] (28) at (2.25, 2) {};
		\node [style=none] (29) at (3, 0) {};
		\node [style=none] (30) at (1.75, -2) {};
		\node [style=none] (31) at (0.5, 2.75) {};
		\node [style=none] (32) at (0.5, 1.35) {};
		\node [style=none] (33) at (0.5, 1) {};
		\node [style=none] (34) at (0.75, 0) {};
		\node [style=none] (35) at (-2.5, -0.5) {};
		\node [style=none] (36) at (-0.25, -2) {};
		\node [style=none] (37) at (3.25, 3.25) {g};
		\node [style=none] (38) at (-2.25, 3.25) {h};
		\node [style=none] (39) at (-2.25, -2.25) {e};
		\node [style=none] (40) at (3.75, -2) {f};
	
		\draw (1.center)
			 to [bend left=60, looseness=0.75] (2.center)
			 to [bend left=60, looseness=0.50] (0.center)
			 to [bend left=60, looseness=0.50] cycle;
		\draw (4.center)
			 to [bend right=60, looseness=0.75] (5.center)
			 to [bend right=60, looseness=0.50] (3.center)
			 to [bend right=60, looseness=0.50] cycle;
		\draw (7.center)
			 to [bend right=60, looseness=0.75] (8.center)
			 to [bend right=60, looseness=0.50] (6.center)
			 to [bend right=60, looseness=0.50] cycle;
		\draw (10.center)
			 to [bend left=60, looseness=0.75] (11.center)
			 to [bend left=60, looseness=0.50] (9.center)
			 to [bend left=60, looseness=0.50] cycle;
		\draw [style=new edge style 1,fill=cl1] (26.center)
			 to [bend right=90, looseness=0.50] (25.center)
			 to [bend left=270, looseness=0.75] cycle;
		\draw [style=new edge style 1,fill=cl1] (28.center)
			 to [bend right=90, looseness=1.25] (27.center)
			 to [bend right=90, looseness=1.50] cycle;
		\draw [style=new edge style 1,fill=cl1] (30.center)
			 to [bend right=90, looseness=1.25] (29.center)
			 to [bend left=270] cycle;
		\draw [style=new edge style 1,fill=cl1] (32.center)
			 to [bend right=90, looseness=1.25] (31.center)
			 to [bend right=90, looseness=1.50] cycle;
		\draw [style=new edge style 1,fill=cl1] (34.center)
			 to [bend right=90, looseness=1.25] (33.center)
			 to [bend right=90, looseness=1.50] cycle;
		\draw [style=new edge style 1,fill=cl1] (36.center)
			 to [bend right=60] (35.center)
			 to [bend right=90] cycle;
   
	    \node [style=new style 0,scale=0.7] (12) at (0.5, 0.5) {1};
		\node [style=new style 0,scale=0.7] (13) at (-2, -1) {2};
		\node [style=new style 0,scale=0.7] (14) at (-1, -1.75) {3};
		\node [style=new style 0,scale=0.7] (15) at (2, -1.5) {4};
		\node [style=new style 0,scale=0.7] (16) at (2.75, -0.5) {5};
		\node [style=new style 0,scale=0.7] (17) at (0.4, 2.5) {6};
  \node [style=new style 0,scale=0.7] (17') at (0.6, 1.75) {8};
		\node [style=new style 0,scale=0.7] (18) at (2.25, 2.5) {7};
		\node [style=new style 0,scale=0.7] (19) at (-1.25, 2.5) {9};
		\node [style=new style 0,scale=0.7] (20) at (-1.75, 1.5) {10};
		\node [style=none,scale=0.7] (21) at (0.45, 2) {};
  \node [style=none,scale=0.7] (41) at (2, 0.5) {$W_{E_1}$};
		\node [style=none,scale=0.7] (42) at (-1.25, -1) {$W_{E_2}$};
		\node [style=none,scale=0.7] (43) at (-1.25, 2) {$W_{E_6}$};
		\node [style=none,scale=0.7] (44) at (2.75, -1.25) {$W_{E_3}$};
		\node [style=none,scale=0.7] (45) at (3, 2) {$W_{E_4}$};
		\node [style=none,scale=0.7] (46) at (1, 3.25) {$W_{E_5}$};
  	\draw [->] (46.center) to (21);
		\draw [->] (45.north) to (18);
		\draw [->, bend right=60, looseness=1.25] (41.north) to (12);
\end{tikzpicture}

\caption{A hypergraph $H$ with the vertex set $V(H)=\{i\in\mathbb{N}:1\le i\le 10\}$, and hyperedge set $E(H)=\{e,f,g,h\}$, where $e=\{1,2,3\}$, $f=\{1,4,5\}$, $g=\{1,6,7,8\}$, $h=\{1,6,8,9,10\}$. In each grey region, all the vertices are incident with the same set of hyperedges. Each grey region is a unit here. Therefore, each unit corresponds to a subset of hyperedges $E \subseteq E(H)$ such that for all vertices $v$ in the unit, the set of hyperedges incident with $v$ is $E_v(H) = E$. This unit is represented as $W_E$. Here the units are $W_{E_1}=\{1\}$, $W_{E_2}=\{2,3\}$, $W_{E_3}=\{4,5\}$, $W_{E_4}=\{7\}$, $ W_{E_5}=\{6,8\}$, $ W_{E_6}=\{9,10\}$. The corresponding generating sets are $E_1=\{e,f,g,h\}$, $E_2=\{e\}$, $E_3=\{f\}$, $E_4=\{g\}$, $E_5=\{g,h\}$, $E_6=\{h\}$.}
		\label{fig:unit}
	\end{figure}
	\begin{exm}\label{ex-unit-hypg}
		Consider the hypergraph $H$ in \Cref{fig:unit}, with $V(H)=\{i\in\mathbb{N}:1\le i\le 10\}$, and $E(H)=\{e,f,g,h\}$, where $e=\{1,2,3\}$, $f=\{1,4,5\}$, $g=\{1,6,7,8\}$, $h=\{1,6,8,9,10\}$.
		In $H$, the units are $W_{E_1}=\{1\}$, $W_{E_2}=\{2,3\}$, $W_{E_3}=\{4,5\}$, $W_{E_4}=\{7\}$, $ W_{E_5}=\{6,8\}$, $ W_{E_6}=\{9,10\}$. Here the generating sets are $E_1=\{e,f,g,h\}$, $E_2=\{e\}$, $E_3=\{f\}$, $E_4=\{g\}$, $E_5=\{g,h\}$, $E_6=\{h\}$.
	\end{exm}
Now, we present a set of matrices associated with a hypergraph $H$ that captures the effects of symmetry within the units of $H$. We denote the collection of all the stars in a hypergraph $H$ as $ \mathfrak{E}(H)=\{E_i(H):i\in V(H)\}$.
	\begin{df}[Star dependent matrix]\label{star-matrix}Let $H$ be a hypergraph. We call a matrix $M_H=\left(m_{ij}\right)_{i,j\in V(H)}$ as a \emph{star dependent matrix} if there exists two functions $g:\mathfrak{E}(H)\times \mathfrak{E}(H)\to \mathbb{R}$, and $h:\mathfrak{E}(H)\to \mathbb{R}$ such that,
		\[m_{ij}=
		\begin{cases}
			g(E_i(H),E_j(H))& \text{~if~}i\ne j,\\
			h(E_i(H))& \text{~otherwise.~}
		\end{cases}
		\]
  We refer to the functions $h$ and $g$, respectively, as the \emph{diagonal and non-diagonal star-function} of the star-dependent matrix $M_H$, respectively.
	\end{df}
	
	\begin{exm}[Adjacency Matrix]\label{ex:adj}
		\begin{enumerate}[leftmargin=*]
			\item According to \cite{bretto2013hypergraph}, the adjacency matrix of a hypergraph $H$ is a square matrix $A^{(1)}$, whose rows and columns are indexed by $V(H)$ and for $i,j(\ne i)\in V(H)$, the $(i,j)$-th entry is $A^{(1)}_{ij}=|E_i(H)\cap E_j(H)|$, the co-degree of $i$ and $j$, and the diagonal position is $A^{(1)}_{ii}=0$ for all $i\in V(H)$. Since $A^{(1)}_{ij}=|E_i(H)\cap E_j(H)|$ depends only on $E_i(H)$, $E_j(H)$  and $A^{(1)}_{ii}=0$ is a constant, $A^{(1)}$ is a star-dependent matrix.
			
			\item Similarly, the adjacency matrix $A^{(2)}=\left(A^{(2)}_{ij}\right)_{i,j\in V(H)}$ of $H$, introduced in \cite{Banerjee-2021-hgmat}, defined by 
			$$A^{(2)}_{ij}:=
			\begin{cases}
				\sum\limits_{e\in E_i(H)\cap E_j(H)}\frac{1}{|e|-1}&\text{if~} i\neq j\\
				0&\text{otherwise}
			\end{cases}$$
			is also a star-dependent matrix. Here, $g(E_i(H),E_j(H))=\sum\limits_{e\in E_i(H)\cap E_j(H)}\frac{1}{|e|-1}$, and $h(E_i(H))=0$ for all $i,j\in V(H)$.
			\item The normalized adjacency matrix $\mathcal{A}=\left(\mathcal{A}_{ij}\right)_{i,j\in V(H)}$ of $H$, defined by 
			$$\mathcal{A}_{ij}:=
			\begin{cases}
				\frac{1}{|E_i(H)|}\sum\limits_{e\in E_i(H)\cap E_j(H)}\frac{1}{|e|-1}&\text{if~} i\neq j\\
				0&\text{otherwise}
			\end{cases}$$
			is a star-dependent matrix.
		\end{enumerate}
	\end{exm}
	\begin{exm}[Laplacian Matrix]\label{ex:lap}
		\begin{enumerate}[leftmargin=*]
			\item
			The Laplacian matrix $L^{(1)}=D^{(1)}-A^{(1)}$ is a star-dependent matrix, where $D^{(1)}$ is the diagonal matrix of order $|V(H)|$, defined by $D^{(1)}_{ii}=(A^{(1)}\chi_{V(H)})(i)$. This variation of hypergraph Laplacian $L^{(1)}$ is used in \cite{rodriguez2003laplacian}. Here  $g(E_i(H),E_j(H))=-|E_i(H)\cap E_j(H)|$, and $h(E_i(H))=|E_i(H)|$ for all $i,j\in V(H)$.
			\item Another variation of Laplacian $L^{(2)}=D^{(2)}-A^{(2)}$, introduced in \cite{Banerjee-2021-hgmat}, is also  a star-dependent matrix. Here, $g(E_i(H),E_j(H))=-\sum\limits_{e\in E_i(H)\cap E_j(H)}\frac{1}{|e|-1}$, and $h(E_i(H))=\sum\limits_{e\in E_i(H)}\frac{1}{|e|-1}$ for all $i,j\in V(H)$.
			\item The normalized Laplacian $\Delta=I_{|V(H)|}-A^{(2)}$ introduced in \cite{Banerjee-2021-hgmat}, is a star-dependent matrix. Here $I_{|V(H)|}$ is the identity matrix of order $|V(H)|$.
		\end{enumerate}
	\end{exm}
	\begin{thm}\label{unit-finer}
		If $H$ is a hypergraph and $M_H$ is a star-dependent matrix, then $\mathfrak{R}_s(H)$ is finer than $\mathfrak{R}_{M_H}$.
	\end{thm}
	\begin{proof}
		For $i,j\in V(H)$, assuming $(i,j)\in \mathfrak{R}_s(H)$, it is enough to show $(i,j)\in \mathfrak{R}_{M_H} $. Since $(i,j)\in \mathfrak{R}_s(H)$ we have $E_i(H)=E_j(H)$. Since $M_H=\left(m_{ij}\right)_{ij\in V(H)}$ is a star-dependent matrix, $E_i(H)=E_j(H)$ leads us to $m_{ii}=m_{jj}$, $m_{ij}=m_{ji}$, $m_{ki}=m_{kj}$, and $m_{ik}=m_{jk}$ for all $k\in V(H)$. This completes the proof.
	\end{proof}
Suppose that $(i,j)\in\mathfrak{R}_s(H)$. In that case $E_i(H)=E=E_j(H)$. Thus,
for a star-dependent matrix $M_H$, we have $m_{ii}=h(E)=m_{jj}$, $m_{ij}=g(E,E)=m_{ji}$, and for any $k\in V(H)\setminus\{i,j\}$, we have $ m_{ik}=g(E,E_k(H))=m_{jk}$, $m_{ki}=g(E_k(H),E)=m_{kj}$. Thus, $(i,j)\in \mathfrak{R}_{M_H}$.
 That is
 for a star-dependent matrix $M_H$,
 $\mathfrak{R}_s(H)$ is finer than $\mathfrak{R}_{M_H}$, for all $(i,j)\in \mathfrak{R}_s(H)$, $\psi_{ij}$ and $M_H$ commutes, that is $\psi_{ij}(M_H x)= M_H(\psi_{ij}(x))$.
\begin{exm}\label{finer-ex}
    \begin{enumerate}[leftmargin=*]
        \item In \Cref{ex-finer-1}, we have already provided a matrix $M_H$ such that $\mathfrak{R}_s(H)$ is finer than $\mathfrak{R}_{M_H}$. In that case for two distinct vertices $i,j$, the $(i,j)$-th entry of the matrix $m_{ij} = \sum\limits_{e \in E_i(H) \cap E_j(H)} \frac{1}{|e|} $, and the diagonal entries of the matrices are $0$. That is, $M_H$ is star-dependent. Consequently, by \Cref{unit-finer},  $\mathfrak{R}_s(H)$ is finer than $\mathfrak{R}_{M_H}$, and $\psi_{ij}(M_H x)= M_H(\psi_{ij}(x))$ for all $(i,j)\in \mathfrak{R}_s(H)$. 
      \item 
      If we set \(M_H = A^{(1)}\), then for the hypergraph $H$ described in \Cref{ex-unit-hypg} (see \Cref{fig:unit}), the matrix $A^{(1)}$ is given below.
      
      \begin{figure}[ht]
          \centering
          \begin{tikzpicture}[scale=0.5]
              \draw[] (0,0) node {$A^{(1)}$=
      \begin{tabular}{c|c||c|cc|cc|c|cc|cc||}
         &  &$W_{E_1}$&$W_{E_2}$& &$W_{E_3}$& &$W_{E_4}$&$W_{E_5}$&&$W_{E_6}$&\\\hline
         &    &$1$&$2$&$3$&$4$&$5$&$7$&$6$&$8$&$9$&$10$\\ \hline\hline
        $W_{E_1}$ &$1$    &$0$&$1$&$1$&$1$&$1$&$1$&$2$&$2$&$1$&$1$\\ \hline
          $W_{E_2}$ &  $2$  &$1$&$0$&$1$&$0$&$0$&$0$&$0$&$0$&$0$&$0$\\
            &  $3$  &$1$&$1$&$0$&$0$&$0$&$0$&$0$&$0$&$0$&$0$\\\hline
           $W_{E_3}$   &$4$    &$1$&$0$&$0$&$0$&$1$&$0$&$0$&$0$&$0$&$0$\\
               &$5$    &$1$&$0$&$0$&$1$&$0$&$0$&$0$&$0$&$0$&$0$\\\hline
               $W_{E_4}$  & $7$   &$1$&$0$&$0$&$0$&$0$&$0$&$1$&$1$&$0$&$0$\\ \hline
             $W_{E_5}$    &  $6$  &$2$&$0$&$0$&$0$&$0$&$1$&$0$&$2$&$1$&$1$\\ 
                 &  $8$  &$2$&$0$&$0$&$0$&$0$&$1$&$2$&$0$&$1$&$1$\\ \hline
                $W_{E_6}$   & $9$   &$1$&$0$&$0$&$0$&$0$&$0$&$1$&$1$&$0$&$1$\\ 
                   & $10$   &$1$&$0$&$0$&$0$&$0$&$0$&$1$&$1$&$1$&$0$\\ \hline\hline
                  \end{tabular}
   };
          \end{tikzpicture}
      \end{figure}
          For all \((i,j) \in \mathfrak{R}_s(H)\), we have \(E_i(H) = E_j(H)\). This implies for all $k\notin\{i,j\}$, \(E_i(H) \cap E_k(H) = E_j(H) \cap E_k(H)\), meaning \((i,j) \in \mathfrak{R}_{A^{(1)}}\). Thus, \(\mathfrak{R}_s(H)\) is finer than \(\mathfrak{R}_{A^{(1)}}\) and \(\psi_{ij}(A^{(1)} x) = A^{(1)}(\psi_{ij}(x))\) for all \((i,j) \in \mathfrak{R}_s(H)\). Similarly, if we choose any star-dependent matrix described in \Cref{ex:adj} and \Cref{ex:lap} as \(M_H\), then we obtain a similar result.

  \end{enumerate} 
\end{exm}
 
{For any star-dependent matrix $M_H$, \Cref{eigen-lemma-1} and the proof of \Cref{unit-finer} show that if $W_E$ is a unit in $H$, then for all $i,j(\ne i)\in W_E$, $\lambda_{W_E}=M_Hx_{ij}(i)$ is an eigenvalue of $M_H$ and $x_{ij}\in V_{\lambda_{W_E}}(M_H)$. Therefore, $S_{W_E}$ is a subspace of $ V_{\lambda_{W_E}}(M_H)$.}
	Thus, \Cref{finer_class_const} and \ref{unit-finer} lead us to the following result on an invariant subspace of star-dependent matrices.
	\begin{cor}[Unit-invariant subspace]\label{unit-invariant}
		Let $H$ be a hypergraph, $f: V(H)\to (0,\infty)$ be a function, and $M_H$ be an $f$-symmetric, star-dependent matrix indexed by $V(H)$. For any $W\in\mathfrak{U}(H)$ with $|W|\ge 2$, the vector space $C_W^f$ is an invariant subspace of $M_H$.
	\end{cor}
For instance, consider the matrix $M_H=A^{(1)}$ described in \Cref{finer-ex}(2) associated with the hypergraph $H$ described in \Cref{fig:unit}. For the unit $W_{E_5}$, the matrix has eigenvalue $-2$ with eigenvector $x_{68}$ which has values $-1$ and $1$
at vertices $ 6$ and $8$, respectively. For all $v\in V(H)\setminus\{6,8\}$, $x_{68}(v)=0$. The matrix $A^{(1)}$ is a symmetric matrix. That is $A^{(1)}$ is $f$-symmetric where $f(v)=1$ for all $v\in V(H)$. Thus, as \Cref{unit-invariant} indicates,  $C^f_{W_{E_5}}=C_{W_{E_5}}=\{x:V(H)\to\mathbb{R}:x(6)=x(8)\}=\langle x_{68}\rangle^\perp$ is an invariant subspace of $A^{(1)}$. Here $\langle x_{68}\rangle^\perp$ is the perpendicular subspace of the vector space generated by the vector $x_{68}$. Similarly, for the units $W_{E_2}$, $W_{E_3}$ and $W_{E_6}$, we have $C_{W_{E_2}}$, $C_{W_{E_3}}$ and $C_{W_{E_9}}$ are invariant subspaces of $A^{(1)}$.

	In the next result, we show that even if a star-dependent matrix is not $f$-symmetric with respect to any $f$, units provide invariant subspaces.
	\begin{thm}\label{unit_invariant_cluster}
		Let $H$ be a hypergraph and $M_H$ be a star-dependent matrix indexed by $V(H)$. For any $W\in\mathfrak{U}(H)$ with $|W|\ge 2$, the vector space $C_W$ is an invariant subspace of $M_H$.
	\end{thm}
	\begin{proof}
		For any $x\in C_W$, it is enough to show $M_Hx\in C_W$. Since $x\in C_W$, we have $x(i)=c_x$ for all $ i\in W$. 
		
		$M_H=\left(m_{ij}\right)_{i,j\in V(H)}$ is a star-dependent matrix. Thus, if $i,j\in W$, then for any $k\in V(H)\setminus\{i,j\}$,  we have $m_{ik}=m_{jk}$, and $m_{ij}=m_{ji}$. This leads us to 
		\begin{align*}
			(M_Hx)(i)= \sum\limits_{k(\ne j)\in V(H)}m_{ik}x(k)+m_{ij}c_x=\sum\limits_{k(\ne i)\in V(H)}m_{jk}x(k)+m_{ji}c_x= (M_Hx)(j).
		\end{align*}
		Therefore, $M_Hx\in C_W$ for all $x\in C_W$, and this completes the proof.
	\end{proof}

 For example, we have seen in \Cref{finer-ex}, that for the matrix $A^{(1)}$ associated with the hypergraph $H$ given in \Cref{fig:unit}, the equivalence relation  \(\mathfrak{R}_s(H)\) is finer than \(\mathfrak{R}_{A^{(1)}}\). Now, suppose that $x_1:V(H)\to\mathbb{R}$ is a vector such that $x_1(i)=c_1$, a constant for all $i\in W_{E_3}$, and $x_2=A^{(1)}x_1$. For $i,i'\in W_{E_3}$, we have $a_{ii}=a_{i'i'}$, $a_{ii'}=a_{i'i}$, and $a_{ik}=a_{i'k}$  for all $k\in \{i,i'\}$. That leads us to $x_2(i)=x_2(i')$ for all $i,i'\in W_{E_3}$. That is, as \Cref{unit_invariant_cluster} indicates, if $x_1\in C_{W_{E_3}}=\{x:V(H)\to\mathbb{R}:x(4)=x(5)\}$, then $x_2=A^{(1)}x_1\in C_{W_{E_3}}$. A vector $x$ being constant on each unit $W_E$ means the components of the vectors agree with each other in $W_E$, that is $x(i)=x(i')$ for all $i,i'\in W_E$. We refer to this phenomenon as \emph{cluster synchronization} in $W_E$. In \Cref{sec-application}, we further explore the emergence of cluster synchronization in dynamical processes governed by \Cref{random_walk_eqn} and (\ref{diff-eqn}).
The above theorem indicates that units are clusters of vertices such that once they reach a cluster synchronization under the action of a star-dependent matrix, the cluster synchronization is retained in the subsequent steps. That is, for a star-dependent matrix $M_H$ in \Cref{mheqn}, we have $x_n\in C_{W_{E_3}}$ for all $n\ge k$ whenever $x_k\in C_{W_{E_3}}$ for some $k\in \mathbb{N}$.
Similarly, if we consider any star-dependent matrix associated with the hypergraph illustrated in \Cref{fig:unit}, then any unit $W_E$ with cardinality at least $2$ leads to an invariant subspace $C_{W_E}$. That is, for the units $W_{E_2}$, $W_{E_3}$ and $W_{E_6}$, we have $C_{W_{E_2}}=\{x:V(H)\to\mathbb{R}:x(2)=x(3)\}$, $C_{W_{E_3}}=\{x:V(H)\to\mathbb{R}:x(4)=x(5)\}$, $C_{W_{E_6}}=\{x:V(H)\to\mathbb{R}:x(9)=x(10)\}$,  and $C_{W_{E_5}}=\{x:V(H)\to\mathbb{R}:x(6)=x(8)\}$ are invariant subspaces of the star-dependent matrix $M_H$.

In the next result, we prove that given any unit $W_E$ in a hypergraph $H$, and a star-dependent matrix $M_H$ associated with $H$, if for some vector $x:V(H)\to \mathbb{R}$, the vector $y=M_Hx$ is constant on $W_E$, then the vector $x$ is also constant in $W_E$.
\begin{thm}\label{backward shift}
    	Let $H$ be a hypergraph and $M_H$ be a star-dependent matrix indexed by $V(H)$. For any $W_E\in\mathfrak{U}(H)$ with $|W_E|\ge 2$ and $g(E,E)\ne h(E)$ where $g$ and $h$ are, respectively, non-diagonal and diagonal star-functions of $M_H$. If for some vector $x:V(H)\to \mathbb{R}$, the vector $y=M_Hx\in C_W$ then $x\in C_{W_E}$.
\end{thm}
\begin{proof}
Since $y=M_Hx\in C_{W_E}$, for any two $i,i'\in W_E$, we have $y(i)=y(i')$. Since $i,i'$ has the same star we have 
\begin{align*}
   &m_{ii} x(i)+m_{ii'}x(i')+\sum\limits_{j\in V(H)\setminus\{i,i'\}}m_{ij}x(j)=y(i)=y(i')\\&=  m_{i'i'} x(i')+m_{i'i}x(i)+\sum\limits_{j\in V(H)\setminus\{i,i'\}}m_{i'j}x(j).
\end{align*}
Since $M_H$ is star-dependent, and both $i,i'$ have the same star,  we have $m_{ii}=m_{i'i'}$ and $m_{ii'}=m_{i'i}$. Thus, the above equation leads us to
\begin{align*}
    (m_{ii}-m_{ii'})(x(i)-x(i'))=0.
\end{align*}
The condition $g(E,E)\ne h(E)$ leads us to $ m_{ii}\ne m_{ii'}$. Therefore, $x(i)=x(i')$. This completes the proof.
\end{proof}
The condition $g(E, E) \ne h(E)$ is essential for the validity of \Cref{backward shift}. This assumption ensures that the diagonal entries of the matrix differ from the off-diagonal entries. Without this distinction, a matrix with identical diagonal and off-diagonal elements may map a non-synchronized vector to a synchronized one. For example, consider the matrix
$
M =\left[ \begin{smallmatrix}
\frac{1}{2} & \frac{1}{2} \\
\frac{1}{2} & \frac{1}{2}
\end{smallmatrix}\right]
$
and the vector
$
x = \left[\begin{smallmatrix}
\frac{2}{3} \\
\frac{1}{3}
\end{smallmatrix}\right].$
Then,
$
Mx = \left[\begin{smallmatrix}
\frac{1}{2} \\
\frac{1}{2}
\end{smallmatrix}\right],
$
which is a vector with equal components, despite $x$ having unequal components.

\begin{rem}
    To prove \Cref{unit_invariant_cluster}, and \ref{backward shift}, we have only used that if two vertices $i,i'$ having the same star, then for a star dependent matrix  $M_H=(m_{ij})_{i,j\in V(H)}$, we have $m_{ii}=m_{i'i'}, m_{ii'}=m_{i'i}$, and for $k\notin\{i,i'\}$, the entry $m_{ik}=m_{i'k}$. This fact is true not only for units but also for the equivalence class of any equivalence relation that is finer than $\mathfrak{R}_{M_H}$. Therefore, if an equivalence relation $\mathfrak{R}$ is finer than  $\mathfrak{R}_{M_H}$, and $W$ is an $\mathfrak{R}$-equivalence class then similar results holds for any $i,i'\in W$.
\end{rem}
	\subsection{Edge-sum Matrices and their Invariant Subspace}\label{subsec_edgesum}
	In \Cref{star-matrix}, we have introduced star-dependent matrices. Now, we need to introduce a special subclass of star-dependent matrices, where the functions of stars, $g$, and $h $ are sums over the stars.
	\begin{df}[Edge-sum matrix]\label{edge-sum}
		A star-dependent matrix $M_H=\left(m_{ij}\right)_{i,j\in V(H)}$ is called an \emph{edge-sum} matrix 
		if there exists two maps $p_1:E(H)\to\mathbb{R}$, and $p_2:E(H)\to\mathbb{R}$, such that, 
		\[ g(E_i(H),E_j(H))=\sum\limits_{e\in E_i(H)\cap E_j(H)}p_1(e),\text{~and~} h(E_i(H))=\sum\limits_{e\in E_i(H)}p_2(e),\]
  where $g$ and $h$ are, respectively, the non-diagonal and diagonal star-functions of the star-dependent matrix $M_H$.
		For $i,j\in V(H)$, if there exists a bijection $b_{ij}:E_i(H)\to E_j(H)$ such that for any $k\notin W_{E_i(H)}\cup W_{E_j(H)}$,
		\[\sum\limits_{e\in E_i(H)\cap E_k(H)}p_1(e)=\sum\limits_{b_{ij}(e)\in E_j(H)\cap E_k(H)}p_1(b_{ij}(e))\text{~and~} \sum\limits_{e\in E_i(H)}p_2(e)=\sum\limits_{b_{ij}(e)\in E_j(H)}p_2(b_{ij}(e)),\] 
   $\sum\limits_{e\in E_i(H)}p_1(e)=\sum\limits_{b_{ij}(e)\in E_j(H)}p_1(b_{ij}(e))$, and $b_{ij}(E_i(H)\cap E_j(H))=E_i(H)\cap E_j(H)$.
  
  then we say $b_{ij}$ is \emph{compatible with the edge-sum matrix $M_H$}. We refer to $p_1$ and $p_2$ as \emph{non-diagonal} and \emph{diagonal edge-weight}, respectively.
	\end{df}
 For instance, the star dependent matrices $A^{(1)}$, $A^{(2)}$ considered in \Cref{ex:adj}, and the Laplacian $L^{(1)}$, $L^{(2)}$ given in \Cref{ex:lap} are edge-sum Matrices. For $A^{(1)}$, the non diagonal star-function is $p_1(e)=1$ for all $e\in E(H)$, and for $A^{(2)}$ the non-diagonal map $p_1(e)=\frac{1}{|e|-1}$ for all $e\in E(H)$. For both $A^{(1)}$, and $A^{(2)}$, the diagonal maps are zero-maps.
Since $\mathfrak{R}_s(H)$ is an equivalence relation on $V(H)$, given a vertex $i\in V(H)$, there is a unique unit $W_{E_i(H)}$ containing $i$. We will use this fact and a bijection between stars, compatible with an edge-sum matrix, to introduce the following equivalence relation for a hypergraph $H$ and an edge-sum matrix $M_H$.
A bijection between stars, compatible with an edge-sum matrix, leads us to the following equivalence relation for a hypergraph $H$ and an edge-sum matrix $M_H$.
\begin{align*}
   & \mathfrak{R}_2(H,M_H)=\{(i,j)\in V(H)\times V(H):|W_{E_i(H)}|=|W_{E_j(H)}|\\
    &\text{~and there exists a bijection~} b_{ij}\text{~compatible with~}M_H.\}
\end{align*}
Since by the above definition, $(i,i)\in  \mathfrak{R}_2(H,M_H)$ for all $i\in V(H)$, the relation $ \mathfrak{R}_2(H,M_H)$ is reflexive. The ordered pair $(i,j)\in \mathfrak{R}_2(H,M_H)$ implies $|W_{E_j(H)}|=|W_{E_i(H)}|$ and there exists a bijection $b_{ji}=b_{ij}^{-1}$. Consequently, the ordered pair $(j,i)\in \mathfrak{R}_2(H,M_H)$, and the equivalence relation $\mathfrak{R}_2(H,M_H)$ is symmetric. If $(i,j),(j,k)\in \mathfrak{R}_2(H,M_H)$, then (1) $|W_{E_i(H)}|=|W_{E_j(H)}|=|W_{E_k(H)}|$, and (2) 
$b_{ik}=b_{jk}\circ b_{ij}$ is a bijection compatible with $M_H$. Thus, $(i,j),(j,k)\in \mathfrak{R}_2(H,M_H)$ implies $(i,k)\in \mathfrak{R}_2(H,M_H)$, and $\mathfrak{R}_2(H,M_H)$ is transitive. Therefore, $\mathfrak{R}_2(H,M_H)$ is an equivalence relation.

For instance, consider the hypergraph $H$ described in \Cref{exm-MH} (see \Cref{fig:MH}), and choose the matrix $M_H$ as either the matrix described in the same example or any one of the matrices $A^{(1)}$, $A^{(2)}$, $L^{(1)}$, and $L^{(2)}$ described in \cref{ex:adj}, and \Cref{ex:lap}. The $\mathfrak{R}_2(H,M_H)$-equivalence classes are $\{1,2,3\}$, and $\{4\}$. In \Cref{exm-MH}, we have seen these are also the collection of $\mathfrak{R}_{M_H}$-equivalence classes in $H$. That is $\mathfrak{R}_2(H,M_H)$, and $\mathfrak{R}_{M_H}$ coincides for the hypergraph $H$, and for these choices of matrices. However, these relations are not always identical. For instance, with the same choices of $M_H$ consider the hypergraph $H$ described in \cref{ex-unit-hypg} (illustrated in \Cref{fig:unit}), in case of the $\mathfrak{R}_2(H,M_H)$-equivalence class $\{2,3,4,5\}=W_{E_2}\cup W_{E_3}$, it consists of two $\mathfrak{R}_{M_H}$-equivalence classes $\{2,3\}$, and $\{3,4\}$. In the next few results, we describe the interrelation among the three equivalence relations $\mathfrak{R}_s(H)$, $\mathfrak{R}_{M_H}$, and $\mathfrak{R}_2(H,M_H)$, and related invariant subspaces.
Now, we show that $\mathfrak{R}_s(H)$ finer than $\mathfrak{R}_2(H,M_H)$.
	\begin{thm}\label{twin-finer}
		Let $H$ be a hypergraph. If $M_H$ is an edge-sum matrix, then $\mathfrak{R}_s(H)$ finer than $\mathfrak{R}_2(H,M_H)$. 
	\end{thm}
	\begin{proof}
	For $(i,j)\in\mathfrak{R}_s(H)$, we have $W_{E_i(H)}=W_{E_j(H)}$, and the identity map work as $b_{ij}$. Consequently, $b_{ij}$ is compatible with $M_H$. Therefore, $(i,j)\in\mathfrak{R}_s(H)$ implies $(i,j)\in\mathfrak{R}_2(H,M_H)$.
	\end{proof}
 Since $\mathfrak{R}_s(H)$ finer than $\mathfrak{R}_2(H,M_H)$, each $\mathfrak{R}_2(H,M_H)$-equivalence class can be expressed as union of units. In the next results, we show some connections between the relations $\mathfrak{R}_2(H,M_H)$ and  $\mathfrak{R}_{M_H}$.
 \begin{thm}\label{unit-part-singleton}
     Let $H$ be a hypergraph and $M_H$ be an edge-sum matrix. If each of the $\mathfrak{R}_2(H,M_H)$-equivalence class $W$ is such that $W=W_{E_0}$, or  $W=W_{E_0}\cup W_{E_1}\cup\ldots\cup W_{E_k}$ with $|W_{E_i}|=1$ for all $i=0,1,\ldots, k$, then $\mathfrak{R}_2(H,M_H)$ is finer than $\mathfrak{R}_{M_H}$.
 \end{thm}
 \begin{proof}
     If $W=W_{E_0}$, then since $\mathfrak{R}_s(H)$ is finer than $\mathfrak{R}_{M_H}$, there exists a $\mathfrak{R}_{M_H}$-equivalence class $W'$ such that $W\subseteq W'$. If $W=W_{E_0}\cup W_{E_1}\cup\ldots\cup W_{E_k}$ with $|W_{E_i}|=1$ for all $i=0,1,\ldots, k$, then for $i,j\in W$, due to the bijection $b_{ij}$ compatible with $ M_H$, we have $m_{ii}=m_{jj}$, $m_{ij}=m_{ji}$, and $m_{ik}=m_{jk}$, $m_{ki}=m_{kj}$ for all $k\notin\{i,j\}$. Therefore, $(i,j)\in \mathfrak{R}_{M_H}$, and there exists a $\mathfrak{R}_{M_H}$-equivalence class $W'$ such that $W\subseteq W'$. Consequently, $\mathfrak{R}_2(H,M_H)$ is finer than $\mathfrak{R}_{M_H}$.
 \end{proof}
 \Cref{unit-part-singleton}, and \Cref{finer_class_const} provide us with an invariant subspace of edge-sum matrices, which is stated in the corollary below.
\begin{cor}
     Let $H$ be a hypergraph,  $f: V(H)\to (0,\infty)$ be a function, and $M_H$ be an edge-sum matrix. If an $\mathfrak{R}_2(H,M_H)$-equivalence class $W$ with $|W|>1$ is such that either  $W=W_{E_0}$, or  $W=W_{E_0}\cup W_{E_1}\cup\ldots\cup W_{E_k}$ with $|W_{E_i}|=1$ for all $i=0,1,\ldots, k$, then the vector space $C_W^f$ is an invariant subspace of $M_H$.
\end{cor}
 
 Suppose that $(i,j)\in \mathfrak{R}_2(H,M_H)$ for two distinct vertices $i,j$. If $E_i(H)\ne E_j(H)$, then for all $k(\ne i)\in W_{E_i(H)}$, $m_{ik}=\sum\limits_{e\in {E_i(H)}}p_1(e)=\sum\limits_{b_{ij}(e)\in {E_j(H)}}p_1(b_{ij}(e))=m_{jk'}$ for some $k'(\ne j)\in W_{E_j(H)}$. Therefore, if $|W_{E_i(H)}|=|W_{E_j(H)}|$, then $\sum\limits_{k\in W_{E_i(H)}}m_{ik}=m_{ii}+(|W_{E_i(H)}|-1)m_{ik}=m_{jj}+(|W_{E_j(H)}|-1)m_{jk'}=\sum\limits_{k'\in W_{E_j(H)}}m_{jk'}$. Similarly, $m_{ik'}=m_{ij}=m_{kj}$ for all $k(\ne i)\in W_{E_i(H)}$, and $k'(\ne j)\in W_{E_j(H)}$. Therefore,  $|W_{E_i(H)}|=|W_{E_j(H)}|$ leads us to $\sum\limits_{k\in W_{E_i(H)}}m_{jk}=\sum\limits_{k'\in W_{E_j(H)}}m_{ik'}$.
	\begin{thm}\label{edge-sum-cor}
		Let $H$ be a hypergraph, $f: V(H)\to (0,\infty)$ be a function, and $M_H$ be an $f$-symmetric,  edge-sum matrix. For any $\mathfrak{R}_2(H,M_H)$-equivalence class $W$ with $|W|\ge 2$, the vector space $C_W^f$ is an invariant subspace of $M_H$.
	\end{thm}
 \begin{proof}
     Since $\mathfrak{R}_s(H)$ is finer than $\mathfrak{R}_2(H,M_H)$, there exists a collection of units $W_{E_0}, W_{E_1},\ldots, W_{E_k}$ such that $W=W_{E_0}\cup W_{E_1}\cup\ldots\cup W_{E_k}$. Let $y_r=\chi_{W_{E_r}}-\chi_{W_{E_0}}$ for all $r=1,2,\ldots,k$. 

   For any $i\notin W_{E_r}\cup W_{E_0}$, the $i$-th component of $M_Hy_r$ is $M_Hy_r(i)=\sum\limits_{k\in W_{E_r}}m_{ik}-\sum\limits_{k'\in W_{E_0}}m_{ik'}$. Since $W_{E_r}$, and $W_{E_0}$ are two units inside the same $\mathfrak{R}_2(H,M_H)$-equivalence class, for each $k\in W_{E_r}$, and $k'\in W_{E_0}$, there exists a bijection $b_{kk'}$ compatible with $M_H$ such that $m_{ik}=\sum\limits_{e\in E_i(H)\cap E_k(H)}p_1(e)=\sum\limits_{b_{kk'}(e)\in E_i(H)\cap E_{k'}(H)}p_1(e)=m_{ik'}$. Therefore, $M_Hy_r(i)=0$.
   
     For any $i\in W_{E_r}$, the $i$-th component of $M_Hy_r$ is $M_Hy_r(i)=\sum\limits_{k\in W_{E_r}}m_{ik}-\sum\limits_{k\in W_{E_0}}m_{ik}$. Similarly, for any  $j\in W_{E_0}$, the $j$-th component of $M_Hy_r$ is $M_Hy_r(j)=\sum\limits_{k\in W_{E_r}}m_{jk}-\sum\limits_{k\in W_{E_0}}m_{jk}$. Since $|W_{E_r}|=|W_{E_0}|$, we have $\sum\limits_{k\in W_{E_r}}m_{ik}-\sum\limits_{k\in W_{E_0}}m_{ik}=-(\sum\limits_{k\in W_{E_r}}m_{jk}-\sum\limits_{k\in W_{E_0}}m_{jk})$. Therefore, $y_j$ is an eigenvector of $M_H$ for all $j=1,\ldots,k$.
     The vector space \(
S_W=(\bigoplus\limits_{W_{E_i}\subseteq W,|W_{E_i}|>1}S_{W_{E_i}}) \oplus \langle\{y_1,y_2,\ldots,y_k\}\rangle
\), where for any set of vector $S$, the vector space generated by $S$ is denoted as $\langle S\rangle$. Since the edge-sum matrix $M_H$ is star dependent, if $|W_{E_i}|>1$, then by \Cref{eigen-lemma-1}, and  \Cref{unit-finer}, $S_{W_{E_i}}$ is a collection of eigenvectors of $M_H$, and each $y_i$ is an eigenvector of $M_H$, the vector space $S_W$ is an invariant subspace of $M_H$. Therefore, by \Cref{perp_lem}, $C_W^f$ is an invariant subspace of $M_H$.
 \end{proof}
 For example, in the hypergraph $H$ illustrated in \Cref{fig:MH}, two $\mathfrak{R}_2(H,M_H)$-equivalence classes are $W=\{1,2,3\}$ and $W'=\{4\}$. The normalized adjacency matrix $\mathcal{A}$ described in \Cref{ex:adj}(3), associated with $H$ is 
 
\begin{table}[H]
    \centering$\mathcal{A}=$
    \begin{tabular}{c||ccc|c|}
        & $1$ & $2$ & $3$ & $4$ \\ \hline\hline
        $1$ & $0$ & \(\frac{1}{4}\) & \(\frac{1}{4}\) & \(\frac{1}{2}\) \\
        $2$ & \(\frac{1}{4}\) & $0$ & \(\frac{1}{4}\) & \(\frac{1}{2}\) \\
        $3$ & \(\frac{1}{4}\) & \(\frac{1}{4}\) & $0$ & \(\frac{1}{2}\) \\ \hline
        $4$ & \(\frac{1}{3}\) & \(\frac{1}{3}\) & \(\frac{1}{3}\) & $0$\\\hline
    \end{tabular}.
\end{table}

This matrix is $f$-symmetric, where $f(i)=|E_i(H)|$ for all $i\in V(H)$.
 As \Cref{edge-sum-cor} suggests, $C_W^f=\{x:V(H)\to\mathbb{R}:x(1)|E_1(H)|=x(2)|E_2(H)|=x(3)|E_3(H)|\}=\{x:V(H)\to\mathbb{R}:x(1)=x(2)=x(3)\}$ is an invariant subspace of $\mathcal{A}$.
 
 {Note that if a star-dependent matrix $M_H$ is not an edge-sum matrix but  satisfies $m_{ik}=m_{jk}$, and $m_{ki}=m_{kj}$ for all $i,j\in W(\subseteq V(H))$, then there is an invariant subspace of $M_H$ associated with $W$. We show such an instance in \Cref{twin-unit-random}}.
 
 {
  \Cref{twin-finer} shows that the equivalence relation $\mathfrak{R}_s(H)$ is finer than the equivalence relation $ \mathfrak{R}_2(H,M_H)$ for any edge-sum matrix $M_H$ related to $H$. Thus, each $\mathfrak{R}_2(H,M_H)$-equivalence class is a union of units. Now, if two units are in the same $\mathfrak{R}_2(H,M_H)$-equivalence class, there should be a structural symmetry between the pair of units. The following notion shows a possible structural symmetry between a pair of units belonging to a $\mathfrak{R}_2(H,M_H)$-equivalence class.
 }
	\par{\textbf{Twin units:}} Let $H$ be a hypergraph. Two units $W_{E},W_{F}\in \mathfrak{U}(H)$, are referred to as \emph{twin units} \cite{unit} if there exists a bijection $f_{EF}:E\to F$ defined by
	\[f_{EF}(e)=(e\setminus W_{E})\cup W_{F}.\] We refer to $f_{EF}$ as the \emph{canonical bijection} of the twin units  $W_{E}$, and $W_{F}$.
	\begin{exm}
		In the hypergraph we have considered in \Cref{ex-unit-hypg}(see \Cref{fig:unit}), $W_{E_2}$, and $W_{E_3}$ are a pair of twin units.  
	\end{exm}
	In \Cref{ex:adj}(2), we have considered the adjacency matrix $A^{(2)}$, which is an edge-sum matrix with non-diagonal edge-weight $p_1:E(H)\to\mathbb{R}$ is defined by $p_1(e)=\frac{1}{|e|-1}$ for all $e\in E(H)$, and the diagonal weight $p_2:E(H)\to \mathbb{R}$ defined by $p_2(e)=0$ for all $e\in E(H)$. Let $W_{E_r}$, and $W_{E_s}$ be a pair of twin units in $H$ with $|W_{E_r}|=|W_{E_s}| $. For two distinct $i,j\in W_{E_r}\cup W_{E_s}$, if both $i,j\in W_{E_r} $ or $i,j\in W_{E_s} $ then $b_{ij}$ is the identity map. If $i\in W_{E_r}$, and $j\in W_{E_s}$ then we set $b_{ij}=f_{E_rE_s}$, the canonical bijection of the pair of twin units  $W_{E_r}$, and $W_{E_s}$. Thus, $(i,j)\in \mathfrak{R}_2(H,A^{(2)})$. Therefore, there exists an $\mathfrak{R}_2(H,A^{(2)})$-equivalence class $W$, such that  $W_{E_r}\cup W_{E_s}\subseteq W$. For a positive-valued function $f:V(H)\to (0,\infty)$, by \Cref{edge-sum-cor}, $C^f_W$ is an invariant subspace of $A^{(2)}$. Similarly, for the other matrices considered in \Cref{ex:adj} and \Cref{ex:lap}, we have similar invariant subspaces. Now, the union of twin units is a set of vertices with similar stars.
 Now, we show that the union of equipotent twin units leads us to some invariant subspaces because of structural symmetry. 
  
  We refer to an edge-sum matrix $M_H$ as \emph{edge-cardinality} matrix if the weights $p_1$ and $p_2$  only depend on the cardinality of the hyperedges. That is, there exists two functions $q_1:\mathbb{N}\to\mathbb{R}$, and $q_2:\mathbb{N}\to\mathbb{R}$ such that $p_1(e)=q_1(|e|)$, and $p_2(e)=q_2(|e|)$ for all $e\in E(H)$. For example, the matrices $A^{(1)},A^{(2)}$, and $L^{(1)},L^{(2)}$ described in \Cref{ex:adj}, and \Cref{ex:lap} are edge-cardinality matrices.
	\begin{thm}\label{eigen-lem-2}
		Let $H$ be a hypergraph, and $M_H=(m_{ij})_{i,j\in V(H)}$ be an edge-cardinality matrix associated with $H$. If $W_{E}$, and $W_{F}$ are a pair of distinct twin units such that $|W_{E}|=|W_{F}|$, then there exists an eigenvalue $\lambda$ of $M_H$ with an eigenvector $\chi_{W_{E}}-\chi_{W_{F}}$. 
	\end{thm}
	\begin{proof}
		For any $r\in W_{E}$,
		\[\sum\limits_{k\in W_E}m_{rk}=(|W_{E}|-1)\sum\limits_{e\in E}q_1(|e|)+\sum\limits_{e\in E}q_2(|e|)=\mu_{W_E}
        .\]
		Similarly, we have the constant $\mu_{W_F}$ associated with $W_{F}$. If two units $W_E$, and $W_F$ are distinct, that is $W_E\ne W_F$, then either $W_E\setminus W_F\ne \emptyset$ or $W_F\setminus W_E\ne \emptyset$. Without loss of generality let $W_E\setminus W_F\ne \emptyset$, and $i\in W_E\setminus W_F$. Thus, $E=E_i(H)\ne F$. Therefore, $W_E\ne W_F$ implies $W_E\cap W_F=\emptyset$, otherwise if $W_E\cap W_F\ne \emptyset$, and $k\in W_E\cap W_F$, then $E=E_k(H)=F$, a contradiction. Given that $|W_{E}|=|W_{F}|$, the two units have no common elements, and $\mathfrak{f}_{EF}(e)=(e\setminus W_E)\cup W_F$, it follows that $|e|=|\mathfrak{f}_{EF}(e)|$.
		Thus, 
		\begin{align*}
			\mu_{W_E}&=(|W_{E}|-1)\sum\limits_{e\in E}q_1(|e|)+\sum\limits_{e\in E}q_2(|e|)\\
			&=(|W_{F}|-1)\sum\limits_{e\in E}q_1(|\mathfrak{f}_{EF}(e)|)+\sum\limits_{e\in E}q_2(|\mathfrak{f}_{EF}(e)|)=\mu_{W_F}.
		\end{align*}
		Suppose that $y_{EF}=\chi_{W_E}-\chi_{W_F}$. We consider the following cases to compute $ (M_H y_{EF})$.
		\begin{enumerate}[leftmargin=*]
			\item For $r\in W_E$, if $k\in W_F$, then $E_r(H)\cap E_k(H)=\emptyset$, and thus, $m_{rk}=\sum\limits_{E_r(H)\cap E_k(H)}p_1(e)=0$. Since $y_{EF}(k)=0$ for all $k\in V(H)\setminus(W_E\cup W_F)$ 
			\begin{align*}
				(M_H y_{EF})(r)= \sum\limits_{k\in W_E}m_{rk}=\mu_{W_E}.
			\end{align*}
			\item Similarly, for $r\in W_F$, $ (M_H y_{EF})(r)=-\mu_{W_F}$.
			\item 
			Since $|W_E|=|W_F|$, 
   If $r\in V(H)\setminus (W_E\cup W_F)$, then for any $k\in W_E,s\in W_F$,
			\begin{align*}
				m_{rk}=\sum\limits_{e\in  E_r(H)\cap E}q_1(|e|)=\sum\limits_{\mathfrak{f}_{EF}(e)\in  E_r(H)\cap F}q_1(|\mathfrak{f}_{EF}(e)|)=m_{rs}.
			\end{align*}
			Therefore, for $r\in V(H)\setminus (W_E\cup W_F)$, 
			\begin{align*}
				(M_H y_{EF})(r)= \sum\limits_{k\in W_E}m_{rk}-\sum\limits_{s\in W_F}m_{rs}=0.
			\end{align*}
		\end{enumerate}
		Thus,
		\begin{align}
			(M_Hy_{EF})(r)=
			\begin{cases}
				\phantom{-}\mu_{W_E}&\text{~if~}r\in W_E\\
				-\mu_{W_F}&\text{~if~}r\in W_F\\
				\phantom{-}0&\text{~otherwise.}
			\end{cases}
		\end{align}
		Since $\mu_{W_E}=\mu_{W_F}$, the result follows.
	\end{proof}
	For example, in the hypergraph $H$ illustrated in \Cref{fig:unit}, $W_{E_2}$, and $W_{E_3}$ are a pair of twin units. The matrix $A^{(1)}$ associated with this hypergraph $H$ has an eigenvalue $1$ with eigenvector $$\chi_{W_{E_2}}-\chi_{W_{E_3}}=\left[\begin{smallmatrix}
	    \phantom{-} 0\\ -1\\ -1\\ \phantom{-} 1\\  \phantom{-}1\\ \phantom{-} 0\\ \phantom{-} 0\\ \phantom{-} 0\\  \phantom{-}0\\  \phantom{-}0
	\end{smallmatrix}\right]$$ due to this pair of twin units. As indicated in \Cref{eigen-lem-2}, the eigenvector due to this pair of twin units will be the same for any other edge-cardinality matrix. Therefore, the matrix $A^{(2)}$ has the eigenvalue $\frac{1}{2}$ with the same eigenvector due to this pair of twin units.
	Now, we show that the union of a pair of equipotent twin units provides us with an invariant subspace of any edge-cardinality matrix associated with a hypergraph.
	\begin{thm}\label{twin-unit-invariant}
		Let $H$ be a hypergraph, $f: V(H)\to (0,\infty)$ be a function, and $M_H$ be an $f$-symmetric,  edge-cardinality matrix. If $W_E$, and $W_F$ are a pair of distinct twin units with $|W_E|=|W_F|$, then $C^f_W$ is an invariant subspace of $M_H$, where $W=W_E\cup W_F$.
	\end{thm}
	\begin{proof}
		Since each edge cardinality matrix is a star dependent matrix, by \Cref{unit-finer}, and \Cref{eigen-lemma-1}, there exists $\lambda_{W_E},\lambda_{W_F}\in \sigma_{M_H}$ such that, $S_{W_E}$ and $S_{W_F}$ are subspaces of $ V_{\lambda_{W_E}}(M_H)$, and $ V_{\lambda_{W_F}}(M_H)$ respectively. Thus, both $S_{W_E}$, and $S_{W_F}$ are invariant subspaces of $M_H$. By \Cref{eigen-lem-2}, the vector space $\langle\chi_{W_E}-\chi_{W_F}\rangle$, spanned by the vector $\chi_{W_E}-\chi_{W_F}$ is also an invariant subspace of $M_H$.
  $M_H$ is $f$-symmetric. Thus, $C^f_W=(S_{W_E}\oplus S_{W_F}\oplus \langle \chi_{W_E}-\chi_{W_F}\rangle )^{(\perp,f)}$ is an invariant subspace of $M_H$, where $W=W_E\cup W_F$. 
	\end{proof}
Since $A^{(1)}$, and $A^{(2)}$ are symmetric matrices, we can assume $f$ is such that $f(i)=1$ for all $i\in V(H)$, and by \Cref{twin-unit-invariant}, in the hypergraph $H$ illustrated in \Cref{fig:unit} due to the pair of twin units $W_{E_2}$, and $W_{E_3}$, the vector space $C_{W_{E_2}\cup W_{E_3}}=\{x:V(H)\to \mathbb{R}:x(2)=x(3)=x(4)=x(5)\}$ is an invariant subspace for both the matrices $A^{(1)}$, and $A^{(2)}$ associated with $H$.

 {We have seen that a pair of twin units is a pair of units in the same $\mathfrak{R}_2(H, M_H)$-equivalence class for some specific matrices $M_H$ related to $H$. Now, we show an $\mathfrak{R}_2(H,M_H)$-equivalence class is not necessarily a union of twin units, and there may be other structural symmetries as well, and by \Cref{edge-sum-cor}, these structural symmetries also corresponds to invariant subspaces of edge-cardinality matrices.

 \begin{exm}
 \begin{figure}[ht]
     \centering
     \begin{tikzpicture}[scale=0.8]
	
		\node [style=none] (0) at (0.75, 0.5) {};
		\node [style=none] (1) at (-1, 2.25) {};
		\node [style=none] (2) at (-5, -1) {};
		\node [style=none] (3) at (-3.75, -2.75) {};
		\node [style=none] (4) at (1.5, -2.5) {};
		\node [style=none] (5) at (3, -0.5) {};
		\node [style=none] (6) at (-2, 2.5) {};
		\node [style=none] (7) at (-2.5, 0) {};
		\node [style=none] (8) at (2.25, -2) {};
		\node [style=none] (9) at (2.25, 0) {};
		\node [style=none] (10) at (-4, -0.25) {};
		\node [style=none] (11) at (-4, -2.25) {};
		\node [style=new style 0,scale=0.6] (12) at (-1, -0.5) {$1$};
		\node [style=new style 0,scale=0.6] (13) at (-0.25, 1.25) {$3$};
		\node [style=new style 0,scale=0.6] (14) at (-1.25, 1.5) {$2$};
		\node [style=new style 0,scale=0.6] (15) at (-3.5, -0.75) {$6$};
		\node [style=new style 0,scale=0.6] (16) at (-3.5, -1.5) {$7$};
		\node [style=new style 0,scale=0.6] (17) at (1.25, -0.5) {$4$};
		\node [style=new style 0,scale=0.6] (18) at (1.25, -1.25) {$5$};
		\node [style=none] (19) at (-1.5, -0.5) {};
		\node [style=none] (20) at (-0.25, -0.25) {};
		\node [style=none] (21) at (-1.75, 1.5) {};
		\node [style=none] (22) at (0.25, 1) {};
		\node [style=none] (23) at (1.25, 0) {};
		\node [style=none] (24) at (1.25, -1.75) {};
		\node [style=none] (25) at (-3.25, -0.25) {};
		\node [style=none] (26) at (-3.25, -2) {};
		\node [style=none,scale=0.7] (27) at (-1.87, -0.2) {$W_{E_1}$};
		\node [style=none,scale=0.7] (28) at (-2.25, -1.5) {$W_{E_4}$};
		\node [style=none,scale=0.7] (29) at (2, -1.5) {$W_{E_3}$};
		\node [style=none,scale=0.7] (30) at (-2.1, 1.3) {$W_{E_2}$};
		\node [style=none,scale=0.7] (31) at (1.25, 1.5) {$e$};
		\node [style=none,scale=0.7] (32) at (-0.75, -3) {$f$};
		\node [style=none,scale=0.7] (33) at (-3.75, 1) {$g$};
		\node [style=none] (34) at (-2.75, -3.75) {H};
	
		\draw (0.center)
			 to [bend right, looseness=0.75] (1.center)
			 to [bend right, looseness=0.50] (2.center)
			 to [bend right, looseness=0.75] (3.center)
			 to [bend right, looseness=0.50] cycle;
		\draw (4.center)
			 to [bend right, looseness=0.75] (5.center)
			 to [bend right, looseness=0.50] (6.center)
			 to [bend right, looseness=0.75] (7.center)
			 to [bend right, looseness=0.50] cycle;
		\draw (8.center)
			 to [bend right, looseness=0.75] (9.center)
			 to [bend right, looseness=0.50] (10.center)
			 to [bend right, looseness=0.75] (11.center)
			 to [bend right, looseness=0.50] cycle;
		\draw [style=new edge style 1,fill=gray!30!white] (20.center)
			 to [bend left=75, looseness=1.50] (19.center)
			 to [bend left=75] cycle;
		\draw [style=new edge style 1,fill=gray!30!white] (22.center)
			 to [bend left=75, looseness=0.75] (21.center)
			 to [bend left=75] cycle;
		\draw [style=new edge style 1,fill=gray!30!white] (24.center)
			 to [bend left=75, looseness=1.50] (23.center)
			 to [bend left=75] cycle;
		\draw [style=new edge style 1, bend left=75, looseness=1.50,fill=gray!30!white] (26.center) to (25.center) to [style=new edge style 1, bend left=75] (25.center) to (26.center);
	\node [style=new style 0,scale=0.6] (12) at (-1, -0.5) {$1$};
		\node [style=new style 0,scale=0.6] (13) at (-0.25, 1.25) {$3$};
		\node [style=new style 0,scale=0.6] (14) at (-1.25, 1.5) {$2$};
		\node [style=new style 0,scale=0.6] (15) at (-3.5, -0.75) {$6$};
		\node [style=new style 0,scale=0.6] (16) at (-3.5, -1.5) {$7$};
		\node [style=new style 0,scale=0.6] (17) at (1.25, -0.5) {$4$};
		\node [style=new style 0,scale=0.6] (18) at (1.25, -1.25) {$5$};
  \node [style=none,scale=0.7] (27) at (-1.87, -0.2) {$W_{E_1}$};
		\node [style=none,scale=0.7] (28) at (-2.25, -1.5) {$W_{E_4}$};
		\node [style=none,scale=0.7] (29) at (2, -1.5) {$W_{E_3}$};
		\node [style=none,scale=0.7] (30) at (-2.1, 1.3) {$W_{E_2}$};
		\node [style=none,scale=0.7] (31) at (1.25, 1.5) {$e$};
		\node [style=none,scale=0.7] (32) at (-0.75, -3) {$f$};
		\node [style=none,scale=0.7] (33) at (-3.75, 1) {$g$};
\end{tikzpicture}
\caption{A hypergraph $H$ with units $W_{E_1}=\{1\}$, $W_{E_2}=\{2,3\}$, $W_{E_3}=\{4,5\}$, and $W_{E_4}=\{6,7\}$. This hypergraph
contains some  $\mathfrak{R}_2(H,M_H)$-equivalence class which are not union of twin units where $M_H$ is any one of the edge-cardinality matrices $A^{(1)}$, $A^{(2)}$, $L^{(1)}$, and $L^{(2)}$ described in \Cref{ex:adj}, and \ref{ex:lap}. }
     \label{fig:ex-3}
 \end{figure}
     Let $H$ be a hypergraph with $V(H)=\{n\in\mathbb{N}:1\le n\le7\}$, and $E(H)=\{e,f,g\}$, where $e=\{1,2,3,4,5\}$, $f=\{1,4,5,6,7\}$, $g=\{1,2,3,6,7\}$ (see \Cref{fig:ex-3}). The hypergraph contains $4$ units, $W_{E_1}=\{1\}$, $W_{E_2}=\{2,3\}$, $W_{E_3}=\{4,5\}$, and $W_{E_4}=\{6,7\}$. If $M_H$ is one of the matrices $A^{(1)}$, $A^{(2)}$, $L^{(1)}$, and $L^{(2)}$ described in \Cref{ex:adj}, and \ref{ex:lap}, then the $\mathfrak{R}_2(H,M_H)$-equivalence classes in $H$ are $W_1=\{2,3,4,5,6,7\}$, and $W_2=\{1\}$. None of these are unions of twin units.
 \end{exm}

 } 
	\section{Application}\label{sec-application}
	In the preceding sections, we have seen hypergraph symmetries represented by equivalence relations on the vertex set give rise to invariant subspaces in matrices associated with these hypergraphs. Specifically, we observed that symmetric clusters of vertices, individual units, and unions of equipotent twin units naturally lead to invariant subspaces for this family of matrices. These invariant subspaces, in turn, provide insights into various phenomena associated with dynamical processes governed by these matrices.
Here, we will illustrate how consistent patterns in clusters of vertices, observed during random walks and diffusion on hypergraphs, can be understood through invariant subspaces. First, we adopt some terminology related to dynamical processes from \cite{jost2005dynamical}. 
A \emph{dynamical system} is an ordered pair $(f,X)$, where $f:X\to X$ is a mapping. A \emph{trajectory} with initial condition $x_0\in X$ of the dynamical system is the sequence $\{x_0,x_1,\ldots\}$ in $X$ where $x_{t+1}=f(x_t)$. That is, $x_t=f^t(x_0)$ for all $t=1,2,\ldots$. 
Some dynamical processes, including random walks and diffusion on hypergraphs, can be represented as $(f,\mathbb{R}^{V(H)})$, 
where at any time $t$, the \emph{state of the dynamical process} is represented by a vector \( x_t: V(H) \to \mathbb{R} \). 
Our focus is on linear dynamical processes where $f:\mathbb{R}^{V(H)}\to\mathbb{R}^{V(H)}$ are such that $ f(x) = M_H x$ for all $x\in \mathbb{R}^{V(H)}$, with $ M_H $ is a matrix associated with the hypergraph. That is, the dynamical systems considered here are of the form $(M_H,\mathbb{R}^{V(H)})$, where $M_H$ is a matrix associated with a hypergraph. For example, in  \cite{Banerjee-2021-hgmat}, to describe a random walk, the matrix \( M_H \) is assumed to be the transpose of the normalized adjacency matrix \( \mathcal{A} \) of the hypergraph, as introduced in \Cref{ex:adj}(3). Here, for a random walk on a hypergraph, $M_H$ is chosen from a broader class of star-dependent matrices that includes $\mathcal{A}$. For diffusion, \( M_H \) is one of the Laplacians related to the hypergraph \( H \), as discussed in \Cref{ex:lap}.

We refer to a dynamical system $(f,\mathbb{R}^{V(H)})$ to as \emph{forward invariant} in $A$ if $f(A)\subseteq A$ for some $A\subset \mathbb{R}^{V(H)}$. Similarly, a dynamical system  $(f,\mathbb{R}^{V(H)})$ is \emph{backward invariant} in $A$ if $f^{-1}(A)\subseteq A$.
At any given time $t$, a state  $x_t: V(H) \to \mathbb{R}$ is said to exhibit \emph{synchronization}  if  $x_t(v) = c_t\in\mathbb{R}$ for all vertices  $v \in V(H)$. 
Synchronization does not necessarily occur over the entire vertex set; instead, it may take place within a specific $W \subset V(H)$. Given $W \subseteq V(H)$, we say that the state $x_t$ exhibits \emph{cluster synchronization} within the \emph{cluster synchronizing set} $W$ at time $t$ if $x_t(v) = c_t$ for all $v \in W$. Let the subspace  
$$  
C_W = \{ x \in \mathbb{R}^{V(H)} \mid x(i) = c_x \text{ for all } i \in W \}  
$$  
be such that the dynamical system $(f, \mathbb{R}^{V(H)})$ is forward invariant in $C_W$. That means, if $x_t(i)=c_{x_t}$ for all $i\in W$ for some constant $c_{x_t}$, then there exists a constant $c_{x_{t+1}}$ such that $x_{t+1}(i)=c_{x_{t+1}}$, for all $i\in W$. This implies that if cluster synchronization occurs in cluster synchronizing set $W$ at some time step $t$, then synchronization persists for all subsequent time steps $t' \geq t$. Similarly, if $(f, \mathbb{R}^{V(H)})$ is backward invariant in $C_W$, then the existence of cluster synchronization in $W$ at some time step $t$ ensures that synchronization was present at all previous time steps $t' \leq t$. In this section, we demonstrate that certain dynamical systems exhibit forward or backward invariance within particular subspaces of $\mathbb{R}^{V(H)}$ associated with units and twin units. Consequently, some cluster synchronizations with  units, or specific unions of twin units as cluster synchronizing  sets, are retained in subsequent time steps.
 For a detailed discussion of synchronization, we refer the reader to \cite{sync,chen2012laplacian} and other related literature.

	\subsection{Random walk on hypergraphs} Let $\mathbb{Z}_+$ be the set of non-negative integers.
	A random walk on a hypergraph $H$ is a function $r:\mathbb{Z}_+\to V(H)$ such that for any $t>0$, $r(t)$ depends only on its previous time-state, that is, on $r(t-1)$. For any two $u,v\in V(H)$, the probability of the event $\{r(t)=u\text{~and~}r(t+1)=v\}$ is referred to as the \emph{transition probability} $P_{uv}$ and which is independent of $t$. The matrix $P=\left(P_{uv}\right)_{u,v\in V(H)}$ is called the \emph{probability transition} matrix.
	Suppose that at any time $t\in \mathbb{Z}_+$, the map $x_t:V(H)\to [0,1]$ is the \emph{probability distribution of the random walk at time $t$}, that is, $x_t(u)$ is the probability of the event $\{r(t)=u\}$ for all $u\in V(H)$. Therefore, 
\begin{equation}\label{random_walk_eqn}
        x_{t+1}(u)=\sum\limits_{v\in V(H)}x_t(v)P_{vu},\text{~that is,~} x_{t+1}= P^Tx_t.
    \end{equation}
 That is, for random walks, the probability distribution at time $t$ is the state of the process at time $t$. The random walk on the hypergraph $H$ can be viewed as a dynamical system represented by the pair $(P^T, [0,1]^{V(H)})$, where $[0,1]^{V(H)}$ denotes the space of all functions $x: V(H) \to [0,1]$.
  A random walk is non-lazy if for any $u\in V(H)$, if $r(t)=u$ for some $t\in \mathbb{Z}_+$, then $r(t+1)\ne u$, that is $P_{uu}=0$ for all $u\in V(H)$.
	 A non-lazy random walk with a uniform probability distribution was considered in \cite{Banerjee-2021-hgmat}. In that case, the probability transition matrix is given by $P = \mathcal{A}$. Here, $\mathcal{A}$ denotes the normalized adjacency matrix described in \Cref{ex:adj}(3). Both $\mathcal{A}$ and its transpose are star-dependent matrices.  In the following two results, we consider the transpose of the probability transition matrix, \( P^T \), to be a star-dependent matrix. This assumption encompasses cases such as \( P = \mathcal{A} \), as well as other probability distributions where \( P^T \) is star-dependent. For example, the transition matrix \( P = \frac{1}{2}I + \frac{1}{2}\mathcal{A} \), which represents a lazy random walk, satisfies this condition since \( P^T \) is star-dependent. Thus, our assumption naturally includes this scenario. Now we use \Cref{unit_invariant_cluster} to prove that if at a particular time $s$, the values of the probability distribution $x_s$ become equal in all the vertices in a unit, then equality will be retained in the subsequent time steps.
	\begin{thm}\label{unit-random} Let $H$ be a hypergraph, and $(P^T, [0,1]^{V(H)}))$  be a random walk on $H$, where $P^T$ is a star-dependent matrix.  For any unit $W_E\in \mathfrak{U}(H)$, the random walk is forward invariant in the subspace $C_{W_E}$.
	\end{thm}
	\begin{proof}
    It is enough to show that: For any unit $W_E\in \mathfrak{U}(H)$, if $x_s(i)=p_s$, a constant for all $i\in W_E$ for some $s\in \mathbb{Z}_+$, then $x_{s+1}(i)=p_{s+1}$, a constant for all $i\in W_E$ at time $s+1$.
    
    If $x_s(i)=p_s$, a constant for all $i\in W$ for some $s\in \mathbb{Z}_+$, then $x_s\in C_{W_E}$.
		Since $P^T$ is a star dependent matrix  by \Cref{unit_invariant_cluster}, $x_{s+1}\in C_{W_E}$. Thus, the theorem follows.
	\end{proof}
Through successive applications of \Cref{unit-random}, we deduce that if cluster synchronization occurs within a unit at time \( s \), it persists for all subsequent finite times \( t > s \). \Cref{backward shift} ensures that a similar fact is true for a backward shift of time. That is, if at a particular time $s\ge 1$, the values of probability distribution $x_s$ become equal in all the vertices in a unit, then equality has been there at time $s-1$ in the probability distribution $x_{s-1}$.  
 \begin{thm}\label{unit-backward} Let $H$ be a hypergraph, and $(P^T, [0,1]^{V(H)})$  be a random walk on $H$, where $P^T$ is a star-dependent matrix.  For any unit $W_E\in \mathfrak{U}(H)$ with $g(E,E)\ne h(E)$ where $g$ and $h$ are, respectively, non-diagonal and diagonal star-functions of $P^T$, the random walk is backward invariant in the subspace $C_{W_E}$.
 \end{thm}
 \begin{proof}
     Since $P^T$ is star-dependent matrices, the result directly follows from \Cref{backward shift}. That is, for any unit $W_E\in \mathfrak{U}(H)$, if $x_s(i)=p_s$, a constant for all $i\in W_E$ for some $s\ge1$, then $x_{s-1}(i)=p_{s-1}$, a constant for all $i\in W_E$.
 \end{proof}
 \Cref{unit-backward} implies that the matrix $P^T$ does not, in general, project a vector without cluster synchronization onto one with cluster synchronization, provided $g(E,E) \ne h(E)$, where $g$ and $h$ denote the non-diagonal and diagonal star-functions of $P^T$, respectively.
This condition, along with two vertices inside the same unit, ensures the existence of a submatrix 
$
M = \left[\begin{smallmatrix}
a & b \\
b & a
\end{smallmatrix}\right]
$
within $P^T$, with $a \ne b$. Any vector $\left(\begin{smallmatrix} x & y \end{smallmatrix}\right)$ can be expressed as
$
\left(\begin{smallmatrix} x & y \end{smallmatrix}\right) = c \left(\begin{smallmatrix} 1 & -1 \end{smallmatrix}\right) + d \left(\begin{smallmatrix} 1 & 1 \end{smallmatrix}\right)
$
for some scalars $c$ and $d$. Applying $M$ to this vector yields
$
M \left(\begin{smallmatrix} x & y \end{smallmatrix}\right)^T = c(a - b) \left(\begin{smallmatrix} 1 & -1 \end{smallmatrix}\right)^T + d(a + b) \left(\begin{smallmatrix} 1 & 1 \end{smallmatrix}\right)^T.
$
Since $a \ne b$, the component along $\left(\begin{smallmatrix} 1 & -1 \end{smallmatrix}\right)^T$ cannot be eliminated unless $c = 0$. Therefore, if the original vector $\left(\begin{smallmatrix} x & y \end{smallmatrix}\right)^T$ lacks cluster synchronization (i.e., $c \ne 0$), then $M \left(\begin{smallmatrix} x & y \end{smallmatrix}\right)^T$ also lacks synchronization.
 Repeated applications of \Cref{unit-backward} imply that if a cluster synchronization occurs in the unit \( W_E \) at time \( s \), then the same cluster synchronization must have existed at all preceding times \( t < s \). 
 \begin{exm}The following table (\Cref{tab:unit-rand}) displays the probability distributions of a random walk on the hypergraph $H$ discussed in \Cref{ex-unit-hypg} (\Cref{fig:unit}) in four consecutive time steps. Here the probability transition matrix is $P=\mathcal{A}$. As \Cref{unit-random} indicates, since the initial probability distribution (see row 1 in \Cref{tab:unit-rand}) is constant on the units $W_{E_2}=\{2,3\}$, $W_{E_3}=\{4,5\}$, $W_{E_5}=\{6,8\}$ the equality is retained in the subsequent time steps. The initial probability is not constant on the unit $W_{E_6}=\{9,10\}$. Thus, the distributions are not equal in $W_{E_6}$ in subsequent time steps, as suggested by \Cref{backward shift}. Though the distributions are not equal in $W_{E_6}$ in subsequent time steps, their difference decreases in each time step. We explain this phenomenon in \Cref{assymp-syn}.
    \begin{table}[ht]
        \centering
        \begin{tabular}{|c|c c c c c c c c c c|}\hline
        \backslashbox{Time}{Vertex}& 1& 2& 3& 4& 5& 6& 7& 8& 9& 10\\\hline
   1 & 0.0270 &   0.0541 &   0.0541 &   0.0811 &   0.0811 &   0.1081&    0.1351&    0.1081&    0.1622&    0.1892\\\hline
   2 &0.3311&
    0.0304&
    0.0304&
    0.0439&
    0.0439&
    0.1683&
    0.0383&
    0.1683&
    0.0760&
    0.0693\\\hline
   3 &  0.2216&
    0.0566&
    0.0566&
    0.0634&
    0.0634&
    0.1465&
    0.0837&
    0.1465&
    0.0801&
    0.0818\\\hline
   4& 
    0.2738&
    0.0560&
    0.0560&
    0.0594&
    0.0594&
    0.1434&
    0.0673&
    0.1434&
    0.0709&
    0.0705\\\hline
        \end{tabular}
        \caption{The probability distributions for a random walk on the hypergraph \( H \) are shown in \Cref{fig:unit}. Each row in the table represents the probability distribution at discrete time step \( i \). The columns correspond to the vertices of \( H \). The probabilities are approximated up to four decimal places. The initial probability distribution is constant on the units  $W_{E_2}$, $W_{E_3}$, $W_{E_5}$, and the equality is retained in the subsequent time steps. }
        \label{tab:unit-rand}
    \end{table}
 \end{exm}
	Now, we show that a similar result is also true for a union of equipotent twin units. For our next result, we assume the probability transition matrix, $P=\mathcal{A}$, the normalized adjacency matrix described in \Cref{ex:adj}(3).
	\begin{thm}\label{twin-unit-random}
    Let $H$ be a hypergraph, and $(\mathcal{A}^T, [0,1]^{V(H)}))$  be a random walk on $H$. For any pair of twin units $W_E,W_F\in\mathfrak{U}(H)$ with the canonical bijection ${f}_{EF}:E\to F$, if  $|W_E|=|W_F|$ then the random walk is forward invariant in $C_{ W_E\cup W_F}$.
	\end{thm}
	\begin{proof}
    It is enough to show that if $x_s(r)=p_s$ for all $r\in W_E\cup W_F$ for some $s\in \mathbb{Z}_+$, then $x_{s+1}(r)=p_{s+1}$ for all $r\in W_E\cup W_F$.
		Let $W_E,W_F\in\mathfrak{U}(H)$ be a pair of twin units with $|W_E|=|W_F|$ with the canonical bijection ${f}_{EF}:E\to F$. Since $|W_E|=|W_F|$, we have $|{f}_{EF}(e)|=|e|$ for all $e\in E$. For any $i\in W_E$, and $j\in W_F$, we have $\mathcal{A}_{ki}=\frac{1}{|E_k(H)|}\sum\limits_{e\in E_k(H)\cap E_i(H)}\frac{1}{|e|-1}=\frac{1}{|E_k(H)|}\sum\limits_{\mathfrak{f}_{EF}(e)\in E_k(H)\cap E_j(H)}\frac{1}{|e|-1}=\mathcal{A}_{kj}$, where $k\in V(H)\setminus (W_E\cup W_F)$. For any $i,k(\ne i)\in W_E$, and $j, k'(\ne j)\in W_F $ we have $\mathcal{A}_{ki}=\frac{1}{|E|}\sum\limits_{e\in E}\frac{1}{|e|-1}=\frac{1}{|F|}\sum\limits_{\mathfrak{f}_{EF}(e)\in F}\frac{1}{|e|-1}=\mathcal{A}_{k'j} $. Thus, if $x_s(r)=p_s$ for all $r\in W_E\cup W_F$ for some $s\in \mathbb{Z}_+$, then for any pair $i\in W_E$, and $j\in W_F$ we have 
		\begin{align*}
			x_{s+1}(i)&=\sum\limits_{k\in V(H)}x_s(k)\mathcal{A}_{ki}\\
			&=\sum\limits_{k(\ne i)\in W_E}p_s\mathcal{A}_{ki}+\sum\limits_{k\in  V(H)\setminus (W_E\cup W_F)}x_s(k)\mathcal{A}_{ki}\\
			&=\sum\limits_{k'(\ne j)\in W_E}p_s\mathcal{A}_{k'j}+\sum\limits_{k\in  V(H)\setminus (W_E\cup W_F)}x_s(k)\mathcal{A}_{kj}=x_{s+1}(j).
		\end{align*}
		If both $i,j$ are either in $W_E$ or belong to $W_F$, then, since $x_s$ is constant on both the unit, by \Cref{unit-random}, $x_{s+1}(i)=x_{s+1}(j)$. Now, when $i\in W_E$, and $j\in W_F$, we have $x_{s+1}(i)=x_{s+1}(j)$. Therefore, $x_{s+1}$ is constant on $W_E\cup W_F$.
	\end{proof}
    The previous result shows that for any pair of twin units $W_E,W_F\in\mathfrak{U}(H)$ with the canonical bijection ${f}_{EF}:E\to F$, if  $|W_E|=|W_F|$ and there is a cluster synchronization in $W_E\cup W_F$ in probability distribution $x_t$ at time $t$, then the cluster synchronization will be retained in the next time step.
    In the previous result, we have assumed the probability transition matrix, $P=\mathcal{A}$. However, by \Cref{twin-unit-invariant}, if the transpose of the probability transition matrix $P^T$ is an $f$-symmetric, edge-cardinality matrix, then $C^f_{W_E\cup W_F}$ is an invariant subspace of $P^T$. If the function $f:V(H)\to(0,\infty)$ is constant on $W_E\cup W_F$ then $C^f_{W_E\cup W_F}=C_{W_E\cup W_F}$. Therefore, in that case, instead of $P=\mathcal{A}$  we can assume $P^T$ to be an $f$-symmetric, edge-cardinality matrix. For instance, if $f$ is such that $f(i)=|E_i(H)|$ for all $i\in V(H)$, then $f(i)=|E|=|F|$, a constant for all $i\in W_E\cup W_F$.
 \begin{exm}
      Consider a random walk on the hypergraph \(H\) discussed in \Cref{ex-unit-hypg} (\Cref{fig:unit}), starting with a random initial distribution (see the first row of \Cref{tab:rand-twin}). The table shows that in the initial time step, the probability distribution of all the vertices in $W_{E_2}\cup W_{E_3}=\{2,3,4,5\}$ are equal, and as \Cref{twin-unit-random} suggests, this equality is preserved in subsequent time steps. Note that $8\notin W_{E_2}\cup W_{E_3}$ and initial probability at $8$ is also equal with the vertices in $W_{E_2}\cup W_{E_3}$ but this equality is not retained in the subsequent time steps.
      \begin{table}[ht]
          \centering
          \begin{tabular}{|c|c c c c c c c c c c|}
          \hline
 \backslashbox{Time}{Vertex}&1&2&3&4&5&6&7&8&9&10\\\hline         
 1&   0.0303&    0.0606&    0.0606 &   0.0606&    0.0606&    0.1212&    0.1515&    0.0606&    0.1818&    0.2121\\\hline
 2&  0.3232&
    0.0341&
    0.0341&
    0.0341&
    0.0341&
    0.1711&
    0.0328&
    0.1887&
    0.0776&
    0.0701\\\hline
 3&  
    0.2210&
    0.0574&
    0.0574&
    0.0574&
    0.0574&
    0.1501&
    0.0869&
    0.1449&
    0.0827&
    0.0846\\\hline
          \end{tabular}
          \caption{Probability distributions associated with a random walk on the hypergraph $H$ given in \Cref{fig:unit}. In the initial probability distribution (the first row), the value of the distribution is equal in all the vertices that belong to the union of twin units $W_{E_2}\cup W_{E_3}$, and the equality is preserved in the subsequent time steps.}
          \label{tab:rand-twin}
      \end{table}
 \end{exm}
  	
\subsection{Diffusion on hypergraphs}
	We have described some Laplacians related to a hypergraph $H$ in \Cref{ex:lap}. Let $L_H$ be any Laplacian described in that example. If the hypergraph is connected, then the only non-positive eigenvalue of $L_H$ is $0$ with multiplicity $1$. The eigenspace $V_0(L_H)$ corresponding to the $0$ eigenvalue is the $1$-dimensional subspace generated by $\chi_{V(H)}$, the characteristic function of $V(H)$. Thus, the $1$-dimensional subspace generated by $\chi_{V(H)}$ is an invariant subspace of $L_H$. For this reason, the Laplacian is said to have a diffusive influence on the hypergraph $H$. The following difference equation can describe the corresponding diffusion process, which is the dynamical system $(L_H,\mathbb{R}^{V(H)})$.
	
	\begin{equation}\label{diff-eqn}
		x_{t+1}=Lx_t.
	\end{equation}
	It has been observed that if the initial state $x_0$ is of the form $c\chi_{V(H)}+\delta y$ for sufficiently small $\delta>0$, and $y\in \mathbb{R}^{V(H)}$, we have $\lim\limits_{t\to\infty}x_t=c_0\chi_{V(H)}$, that is, the states of all the nodes will reach to a synchronization under the influence of $L_H$. Before reaching a complete synchronization that is $x_t=c\chi_{V(H)}$ for some real $c$, we observed some cluster synchronizations in different subsets of vertices. While exploring the reason for the small clusters of synchronization, we find that units preserve cluster synchronization. In the next result, we show that the dynamical system $(L_H,\mathbb{R}^{V(H)})$ is forward invariant in $C_W$ for any unit $W$.
	\begin{thm}
    Let $H$ be a hypergraph, and $W\in \mathfrak{U}(H)$. The diffusion process  $(L_H,\mathbb{R}^{V(H)})$ is forward invariant in $C_W$.
        
	\end{thm}
	\begin{proof}
		Since $L_H$ is a star-dependent matrix, the result follows from \Cref{unit_invariant_cluster}.
	\end{proof}
  The above result guarantees that if the diffusion process $(L_H,\mathbb{R}^{V(H)})$ achieves cluster synchronization in a unit at any given time, this synchronized state is maintained at all later finite times. 
	A similar result can be proved for a union of equipotent twin units. The proof of the result is similar to the proof of \Cref{twin-unit-random}; thus, we omit the proof.
	\begin{thm}
     Let $H$ be a hypergraph, and $W_E,W_F\in \mathfrak{U}(H)$ are a pair of twin units with $|W_E|=|W_F|$. The diffusion process  $(L_H,\mathbb{R}^{V(H)})$ is forward invariant in $C_{W_1\cup W_2}$.
	\end{thm}

\subsection{Asymptotic cluster synchronization}\label{assymp-syn}
The \Cref{assymp-syn} motivates another notion of synchronization called \emph{asymptotic cluster synchronization}. In asymptotic cluster synchronization, for any pair of vertices $i,j\in V(H)$, the difference between the $i$-th and $j$-th components of the state $x_t$ tends to $0$ as $t\to\infty$. More precisely, the trajectory  $\{x_t : t \in \mathbb{Z}_+\}$ of a dynamical system $(f,X)$ is said to attain \emph{eventual cluster synchronization} or \emph{asymptotic cluster synchronization} with the
cluster synchronizing set $ S $ if $\lim\limits_{t \to \infty} |x_t(i) - x_t(j)| = 0 $ for all $i, j \in S $. Asymptotic cluster synchronization in $ S $ can occur depending on the \emph{initial state} $ x_0 $ and the structure of the hypergraph.
In the subsequent results, we consider the pair $(W, M_H)$, where $W$ is a unit and $M_H$ is a matrix whose corresponding eigenvalue $\lambda_W$ satisfies $|\lambda_W| < 1$. For such a pair, the dynamical system $(M_H, \mathbb{R}^{V(H)})$ exhibits asymptotic cluster synchronization.
\begin{thm}\label{cluster-asymptotic}
    Let $H$ be a hypergraph and $M_H$ is a matrix associated with the hypergraph. If an $\mathcal{R}_{M_H}$-equivalence class $W$ with $|W|>1$ is such that $|m_{ij}-m_{ii}|<1$ then there will be
     an asymptotic cluster synchronization with the cluster synchronizing set $W$ in any trajectory of the dynamical system $(M_H,\mathbb{R}^{V(H)})$.
\end{thm}
\begin{proof}
    In the proof of \Cref{eigen-lemma-1}, we have proved that for distinct $i,j\in W$, $\lambda_W=-(M_Hx_{ij})(i)=m_{ij}-m_{ii}$ is an eigenvalue of $M_H$ and $S_W$ is an subspace of the eigenspace of $\lambda_W$. The condition $|m_{ij}-m_{ii}|<1$ leads to $|\lambda_W|<1$.

    For any trajectory $\{x_t:t=0,1,2,\ldots\}$ of the dynamical system $(M_H,\mathbb{R}^{V(H)})$ with initial state $x_0:V(H)\to \mathbb{R}$, we have $x_0 = y + z $, where \( y \) is an eigenvector corresponding to the eigenvalue \( \lambda_W \) associated with the unit \( W \).  Therefore, \( x_t = \lambda_W^t y + z_t \) where $z_t=M_H z \in C_W$  for all $t=0,1,2,\ldots$. Since $|\lambda_W|<1$ we have $\lim\limits_{t\to\infty}\lambda^t=0$. Since $z_t(i)=z_t(j)$ for all $i,j\in W$, therefore $\lim\limits_{t\to\infty}|x_t(i)-x_t(j)|=0$. This completes the proof.
\end{proof}
\begin{exm}Consider the hypergraph $H$ depicted in \Cref{fig:MH}. Let us examine the dynamical system $(A^{(2)}, \mathbb{R}^{V(H)})$, where $A^{(2)}$ denotes the adjacency matrix defined in \Cref{ex:adj}(2). In this case, we observe that $|m_{ij} - m_{ii}| = 0.5 < 1$ for all $i,j \in \{1,2,3\}$. Therefore, by the criterion provided in \Cref{cluster-asymptotic}, the difference between the states of the system corresponding to the vertices in the set $S = \{1,2,3\}$ diminishes over time (see \Cref{tab:my_label-3}). Hence, the system exhibits asymptotic cluster synchronization with cluster synchronizing set $S$.
\begin{table}[ht]
    \centering
    \begin{tabular}{|c|ccccc|}\hline
    \backslashbox{Vertex}{Time}&1&2&3&4&5\\\hline
      1& 0.250 & 1.625 & 2.937 & 7.531 & 16.484\\
      2&0.500 &1.500 & 3 & 7.500 & 16.500\\ 
      3&0.750 & 1.375 &3.0625 & 7.469 & 16.516\\
      4&1 & 1.500 &4.500& 9 & 22.500\\\hline
    \end{tabular}
    \caption{The above table shows the states of dynamical system $(A^{(2)}, \mathbb{R}^{V(H)})$, where $A^{(2)}$ denotes the adjacency matrix defined in \Cref{ex:adj}(2). The rows of the tables correspond to the vertices. The columns of the table correspond to time-steps of a dynamical system. The first three rows of the table indicate that the differences between the components of the states are decreasing.}
    \label{tab:my_label-3}
\end{table}
    \end{exm}
    The next result directly follows from \Cref{cluster-asymptotic}.
    \begin{cor}
         Let $H$ be a hypergraph and $M_H$ be a star-dependent matrix associated with $H$ such that $g$ and $h$ are, respectively, the non-diagonal and diagonal star-function of the matrix. If $W_E$ is a unit in $H$ with $|W_E|>1$ is such that $|g(E,E)-h(E)|<1$ then there will be
     an asymptotic cluster synchronization with the cluster synchronizing set $W$ in any trajectory of the dynamical system $(M_H,\mathbb{R}^{V(H)})$.
    \end{cor}
    \begin{proof}
       Suppose that $M_H=\left(m_{ij}\right)_{i,j\in V(H)}$. Since $|g(E,E)-h(E)|=|m_{ij}-m_{ii}|$ for any two distinct vertices $i,j\in W_E$, the result follows from \Cref{cluster-asymptotic}.
    \end{proof}
	
    \section*{Acknowledgement} The work of the author, Anirban Banerjee, is supported by the FIST Program of the Department of Science and Technology, Government of India (Project No. SR/FST/MS-II/2019/51). 
    The author, Samiron Parui, expresses gratitude to the National Institute of Science Education and Research for their financial support through the Department of Atomic Energy plan project RIA 4001 (NISER). The authors also extend heartfelt thanks to the anonymous referees for their valuable comments and suggestions, which have greatly enhanced the quality of this work.
	\bibliographystyle{siam}

 \end{document}